%% file: ms.tex
 \newcommand{\classstyle}{0}

\ifnum\classstyle=1
	\documentclass[review,onefignum,onetabnum]{siamonline181217}
\else 
	\documentclass[12pt]{amsart}
\fi

\ifnum\classstyle=2
	\usepackage{wiaspreprint}
\else
    \usepackage[english]{babel}
\fi

\usepackage{epsfig}
\usepackage{amsmath,amsfonts,amssymb,rotating,graphicx,mathtools}

\usepackage{tikz}
\usepackage{pgfplots}
\usepackage{pgfplotstable}
\pgfplotsset{compat=1.7}
\usetikzlibrary{arrows,backgrounds}
\usepackage{comment}
\usepackage{subfigure}
\usepackage{placeins}
\usepackage{enumitem}
\setlist[enumerate]{leftmargin=.5in}
\setlist[itemize]{leftmargin=.5in}
\usepackage[nocompress]{cite}
\usepackage{booktabs}
\usepackage[noexternal]{tikzconditional}
\usepackage{graphicx}

\usepackage[textsize=tiny, disable]{todonotes}

\ifnum\ifnum\classstyle=0 0 \else\ifnum\classstyle=2 0 \else 1\fi\fi = 0 %
  \usepackage{amsthm}
  \theoremstyle{plain}             
  \newtheorem{theorem}{Theorem}[section]
  \newtheorem{assumption}[theorem]{Assumption}
  \newtheorem{remark}[theorem]{Remark}

  \newtheorem{corollary}[theorem]{Corollary}

  \usepackage{algorithm}

\else
  \ifnum\classstyle=1
	\newtheorem{assumption}[theorem]{Assumption}
	\newtheorem{remark}[theorem]{Remark}  
  \fi
\fi

\usepackage{algorithmic}

\input{macros}

\begin{document}

\ifnum\ifnum\classstyle=0 0 \else\ifnum\classstyle=2 0 \else 1\fi\fi = 0 %
  \ifnum\classstyle=2
    \include{title}		  
  \fi
  \ifnum\classstyle=0
    \include{title-main}  
  \fi
  \include{abstract}
  \maketitle
\else
  \include{title-siam}    
  \maketitle
  \input{abstract}

\fi

\input{introduction}
\input{randomDomains}
\input{galerkin}
\input{tensorTrains}
\input{adaptivity}

\input{numericalResults}
\ifnum\classstyle=2
  \bibliographystyle{siamplain}
\else
  \bibliographystyle{plain}
\fi
\bibliography{literature}
\end{document}

%% file: macros.tex
\newcommand{\E}{\mathbb{E}}

\newcommand{\isdef}{\mathrel{\mathrel{\mathop:}=}}
\newcommand{\defis}{\mathrel{=\mathrel{\mathop:}}}

\newcommand{\Cov}{\operatorname{Cov}}

\newcommand{\supp}{\operatorname{supp}}

\newcommand{\diam}{\operatorname{diam}}

\renewcommand{\div}{\operatorname{div}}

\newcommand{\de}{\operatorname{d}\!}

\newcommand{\balpha}{{\boldsymbol{\alpha}}}

\DeclareMathOperator{\argmin}{\operatorname{argmin}}

\newcommand{\on}[1]{\operatorname{#1}}
\newcommand{\bs}[1]{{\boldsymbol#1}}

\newcommand{\bbE}{\mathbb{E}}

\def\mcA{{\mathcal{A}}}

\def\mcF{{\mathcal{F}}}
\def\mcI{{\mathcal{I}}}
\def\mcT{{\mathcal{T}}}

\def\mcR{{\mathcal{R}}}
\def\mcS{{\mathcal{S}}}

\def\mcV{{\mathcal{V}}}

\def\mcX{{\mathcal{X}}}
\def\mcY{{\mathcal{Y}}}

\def\bbR{\mathbb{R}}
\def\bbN{\mathbb{N}}

\newcommand*{\jump}[1]{\ensuremath{ [\![ #1 ]\!] }}
\newcommand*{\bigjump}[1]{\ensuremath{ \left[\!\!\left[ #1 \right]\!\!\right] }}

\DeclareMathOperator{\id}{id}

\newcommand*{\abs}[2][{}]{\ensuremath{#1\lvert#2#1\rvert}}
\newcommand*{\norm}[2][{}]{\ensuremath{#1\lVert#2#1\rVert}}

\newcommand*{\mset}[3][{}]{\ensuremath{#1\{ #2 \,;\: #3#1\}}}

\newcommand*{\dd}{\ensuremath{\mathrm{d}}}
\newcommand*{\dx}[1]{\ensuremath{\,\dd{#1}}}

\DeclareMathOperator{\Solve}{Solve}
\DeclareMathOperator{\Estimate}{Estimate}

\DeclareMathOperator{\Refine}{Refine}

\DeclareMathOperator{\Reconstruct}{Reconstruct}

%% file: title-main.tex
\title[Adaptive random domain TT SGFEM]{An Adaptive Stochastic Galerkin Tensor Train Discretization for Randomly Perturbed Domains}
\date{\today}

\author[M.\ Eigel]{Martin Eigel}
\address{Weierstrass Institute\\Mohrenstrasse 39\\D-10117 Berlin\\Germany}
\email{martin.eigel@wias-berlin.de}

\author[M.\ Marschall]{Manuel Marschall}
\address{Weierstrass Institute\\Mohrenstrasse 39\\D-10117 Berlin\\Germany}
\email{manuel.marschall@wias-berlin.de}

\author[M.\ Multerer]{Michael Multerer}
\address{Universit\`a della Svizzera italiana\\ Institute of Computational Science\\ Via Giuseppe Buffi 136\\ CH-6900 Lugano\\ Switzerland}
\email{michael.multerer@usi.ch}

\subjclass[2010]{%
 35R60, 
 47B80, 
 60H35, 
 65C20, 
 65N12, 
 65N22, 
 65J10
}
\keywords{Partial differential equations with random coefficients, tensor representation, tensor train, uncertainty quantification, stochastic finite element methods, adaptive methods, ALS, low-rank, reduced basis methods}

%% file: abstract.tex
\begin{abstract}
A linear PDE problem for randomly perturbed domains is considered in an adaptive Galerkin framework.
The perturbation of the domain's boundary is described by a vector valued random field depending on a countable number of random variables in an affine way.
The corresponding Karhunen-Lo\`eve expansion is approximated by the pivoted Cholesky decomposition based on a prescribed covariance function.
The examined high-dimensional Galerkin system follows from the domain mapping approach, transferring the randomness from the domain to the diffusion coefficient and the forcing.
In order to make this computationally feasible, the representation makes use of the modern tensor train format for the implicit compression of the problem.
Moreover, an a posteriori error estimator is presented, which allows for the problem-dependent iterative refinement of all discretization parameters and the assessment of the achieved error reduction.
The proposed approach is demonstrated in numerical benchmark problems.
\end{abstract}

%% file: introduction.tex
\section{Introduction}
\label{sec:introduction}

Uncertainties in the data for mathematical models are found naturally when dealing with real-world applications in science and engineering.
Being able to quantify such uncertainties can greatly improve the relevance and reliability of computer simulations and moreover provide valuable insights into statistical properties of quantities of interest (QoI).
This is one of the main motivations for the thriving field of Uncertainty Quantification (UQ).

In the application considered in this work, the computational domain is assumed as randomly perturbed.
This e.g.\ can be an appropriate model to incorporate production tolerances into simulations and extract statistical information about how such uncertainties get transported through the assumed model.
Random domain problems have been examined before, see for instance~\cite{xiu,CNT16,HPS16}.
Often, sampling approaches are used to evaluate QoI as e.g.\ has been investigated with a multilevel quadrature for the the domain mapping method in~\cite{HPS16}.
As an alternative, we propose to employ a stochastic Galerkin FEM (SGFEM) to obtain a functional representation of the stochastic solution on the reference domain, which can then be used to evaluate statistical quantities.
For the discretization, a Legendre polynomial chaos basis and first order FE are chosen.
The expansion of the perturbation vector field in a (finite) countable sequence of random variables gives rise to a high-dimensional coupled algebraic system, which easily becomes intractable to numerical methods or results in very slow convergence.
A way to overcome this problem is to utilize model order reduction techniques.
In this work, we make use of the modern tensor train (TT) format~\cite{oseledets2010tt}, which provides an efficient hierarchical tensor representation and is able to exploit low-rank properties of the problem at hand.
Another important technique to reduce computational cost is the use of an adaptive discretization.
In our case, this is based on a reliable a posteriori error estimator, afforded by the quasi-orthogonal approximation obtained by the SGFEM.
With the described error estimator, an iterative adaptive selection of optimal discretization parameters (steering mesh refinement, anisotropic polynomial chaos and tensor ranks) is possible.

For the Karhunen-Lo\`eve expansion of the random vector field, we employ the pivoted Cholesky decomposition derived in~\cite{HPS,HPS14}.
The random coefficient and right-hand side which arise due to the integral transformation are tackled with a tensor reconstruction method.
All evaluations are carried out in the TT format, which in particular allows for the efficient computation of the error estimator as part of the adaptive algorithm.

The paper is structured as follows: The next section introduces the setting and the required assumptions of the random linear model problem.
In particular, a description of the perturbation vector field and the variable transformation is given, converting the random domain problem to a problem with random coefficient and forcing.
Section~\ref{sec:adaptive galerkin} defines the Galerkin finite element discretization of the random coefficient problem in Legendre chaos polynomials.
Moreover, the framework for residual based a posteriori error estimation is described.
The tensor train format used for the efficient computation of the problem is introduced in Section~\ref{sec:tensor trains}.
Section~\ref{sec:adaptive algorithm} lays out the refinement strategy for the Galerkin method, which is based on the evaluation of a reliable a posteriori error estimate in the tensor representation and an appropriate adaptive algorithm.
Numerical examples are discussed in Section~\ref{sec:numerical examples}.


%% file: randomDomains.tex
\section{Diffusion problems on random domains}
In this section, we formulate the stationary diffusion problem on random domains as introduced in \cite{HPS16}.
Let \((\Omega,\Sigma,\mathbb{P})\) denote a complete and separable probability space with \(\sigma\)-algebra \(\Sigma\) and probability measure \(\mathbb{P}\). 
Here, complete means that \(\Sigma\) contains all \(\mathbb{P}\)-null sets.
Moreover, for a given Banach space $E$, we introduce the Lebesgue-Bochner space $L^p(\Omega, \Sigma, \mathbb{P};\mcX)$, $1\le p\le\infty$, which consists of all equivalence classes of strongly measurable functions $v\colon\Omega\to \mcX$ with bounded norm
\[
  \|v\|_{L^p(\Omega, \Sigma, \mathbb{P};\mcX)}\isdef
  \begin{cases}
    \displaystyle{
      \left(\int_\Omega \|v(\cdot,\omega)\|_{\mcX}^p \,
      \de\mathbb{P}(\omega)\right)^{1/p}
    },
    & p<\infty \\[2ex]
    \displaystyle{
      \operatorname*{ess\,sup}_{\omega\in\Omega} \|v(\cdot,\omega)\|_{\mcX}
    },
    & p=\infty.
  \end{cases}
\]
Note that for $p=2$ and $\mcX$ a separable Hilbert space, $L^p(\Omega, \Sigma, \mathbb{P};\mcX)$ is isomorphic to the tensor product space $\mcX \otimes L^2(\Omega, \Sigma, \mathbb{P})$.
We henceforth neglect the dependence on the \(\sigma\)-algebra to simplify the notation.
For an exposition of Lebesgue-Bochner spaces we refer to \cite{HP57}.

In this article, we are interested in computing quantities of interest of the solution to the elliptic diffusion problem
\begin{equation}\label{eq:modProb}
\begin{aligned}
-\div\big(\nabla u(\omega)\big) &= f&&\text{in }D(\omega),\\
u(\omega) &= 0&&\text{on }\partial D(\omega),
\end{aligned}
\end{equation}
for \(\mathbb{P}\)-almost every \(\omega\in\Omega\). 
Note that, the randomness is carried by the open and bounded Lipschitz domain $D\colon\Omega\to\bbR^d$.
It is also possible to consider non-trivial diffusion coefficients or boundary data, see e.g.\ \cite{GP18} for the treatment
of non-homogenous Dirichlet data and \cite{MM18} for random diffusion coefficients.
However, we emphasize that, in order to derive regularity results that
allow for the data sparse approximation of quantities of interest, the data
have to be analytic functions, cf. \cite{HPS16}.
 
In order to guarantee the well posedness of \eqref{eq:modProb}, 
we assume that all data, i.e.\ the loading \(f\) and a possible non-trivial diffusion coefficient, are defined with respect to the \emph{hold-all domain}
\[
\mathcal{D}\isdef\bigcup_{\omega\in\Omega}D(\omega).
\]

For the modelling of random domains, we employ the concept of random 
vector fields. To that end, we assume that there exists a reference
domain $D_{\on{ref}}\subset\mathbb{R}^d$ for \(d=2,3\) with Lipschitz
continuous boundary \(\partial D_{\on{ref}}\) and a random vector field
\[
{\bs V}\colon D_{\on{ref}}\times\Omega\to\mathbb{R}^d
\]
such that \(D(\omega)={\bs V}(D_{\on{ref}},\omega)\).
In addition, we require that \({\bs V}\) is a \emph{uniform \(C^1\)-diffeomorphism}, 
i.e.\ there exists a constant \(C_{\operatorname{uni}}>1\) such that
\begin{equation}\label{eq:unifVfield}
\|{\bs V}(\omega)\|_{C^1(\overline{D_{\on{ref}}};\mathbb{R}^d)},
\|{\bs V}^{-1}(\omega)\|_{C^1(\overline{D_{\on{ref}}};\mathbb{R}^d)}\leq C_{\operatorname{uni}}
\quad\text{for \(\mathbb{P}\)-a.e.\ }\omega\in\Omega.
\end{equation}

In particular, since \({\bs V}\in L^\infty\big(\Omega;C^1(\overline{D_{\on{ref}}})\big)\subset  
L^2\big(\Omega;C^1(\overline{D_{\on{ref}}})\big)\),
the random vector field \({\bs V}\) exhibits a Karhunen-Lo\`eve expansion of the form
 \begin{equation}\label{eq:randomVectorField}
{\bs V}({\hat x},\omega)=\E[{\bs V}]({\hat x})+
\sum_{k=1}^\infty{\bs V}_k({\hat x}){Y}_k(\omega).
\end{equation}
Herein, the expectation is given in terms of the Bochner integral
\[
\E[{\bs V}]({\hat x})\isdef\int_\Omega {\bs V}({\hat x},\omega)\de\mathbb{P}(\omega).
\]
Note that here and henceforth, we denote \({\hat x}\in D_{\on{ref}}\) as material coordinates, in contrast to spatial coordinates \(x\in D(\omega)\).
In particular, there holds \(x={\bs V}({\hat x},\omega)\) for some \({\hat x}\in D_{\on{ref}}\).
The anisotropy, which is induced by the spatial contributions \(\{{\bs V}_k\}_k\), 
describing the fluctuations around the nominal value $\mathbb{E}[{\bs V}]({\hat x})$,
is encoded by
\begin{equation}\label{eq:GammaK}
\gamma_k\isdef\|{\bs V}_k\|_{W^{1,\infty}(D_{\on{ref}};\mathbb{R}^d)}.
\end{equation}
In our model, we shall also make the following common assumptions.
\begin{assumption}\label{ass:KL}\
    \begin{itemize}[leftmargin=1.2cm]
    \item[(i)] The random variables \(Y=\{Y_k\}_k\) take values in \(\Gamma_1\isdef[-1,1]\).
    \item[(ii)] The random variables \(\{Y_k\}_k\) are independent and identically distributed.
    \item[(iii)] The sequence \(\bs\gamma\isdef\{{\gamma}_k\}_k\) is at least in \(\ell^1(\mathbb{N})\).
\end{itemize}
\end{assumption}

In view of this assumption, the Karhunen-Lo\`eve expansion \eqref{eq:randomVectorField}
can always be computed if the expectation \(\E[{\bs V}]\) and the matrix-valued covariance
function 
\[
\Cov[{\bs V}]({\hat x},{\hat x}')\isdef\int_\Omega
\overline{{\bs V}}({\hat x},\omega)
\overline{{\bs V}}^\intercal({\hat x}',\omega)\de\mathbb{P}(\omega)
\]
are known. Herein,
\[
\overline{{\bs V}}({\hat x},\omega)\isdef {\bs V}({\hat x},\omega)-\E[{\bs V}]({\hat x})
\]
denotes the centered random vector field.
The Karhunen-Lo\`eve expansion is based on the spectral decomposition of the integral operator associated to the covariance function, which
can be computed efficiently by means of the pivoted Cholesky decomposition, if the covariance function is sufficiently smooth, cf.\ \cite{HPS,HPS14}.

By an appropriate reparametrization, we can always guarantee that 
\[
\E[{\bs V}]({\hat x})={\hat x}.
\]
Moreover, if we identify the random variables $\{Y_k\}_k$ by their image \({\bs y}\in \Gamma_\infty\isdef \bigtimes_{m\in\bbN} \Gamma_m = [-1,1]^{\mathbb{N}}\),
we end up with the representation
\begin{equation}\label{eq:ParVfield}
    {\bs V}({\hat x},{\bs y})={\hat x}+
    \sum_{k=1}^\infty{\bs V}_k({\hat x}){y}_k.
\end{equation}
For later reference, we also introduce the push-forward measure $\pi_\infty=\mathbb{P}_{\#}Y$ on $\Gamma_\infty$, which will be assumed as a tensor product measure $\pi_\infty = \bigotimes_{m\in\bbN} \pi_m$, where $\pi_m$ is a probability measure on $\Gamma_m=[-1, 1]$. 

The Jacobian of \({\bs V}\) with respect to the material coordinate $\hat{x}$ is given by
\[
{\bs J}({\hat x},{\bs y})={\bs I}+
\sum_{k=1}^\infty{\bs V}_k'({\hat x}){y}_k.
\]
Introducing the parametric domains \(D({\bs y})\isdef{\bf V}(D_{\on{ref}},{\bs y})\), i.e.\
\[
x={\bs V}({\hat x},{\bs y}),
\]
we may now introduce the model problem transported to the reference domain which reads for every $\bs y\in\Gamma_\infty$:
\begin{equation}\label{eq:parmodProb}
\begin{aligned}
-\div_{{\hat x}}\big({\bs A}({\bs y})\nabla_{{\hat x}}\hat{u}({\bs y})\big) &=\hat{f}({\bs y})&&\text{in }D_{\on{ref}},\\
\hat{u}({\bs y}) &= 0&&\text{on }\partial D_{\on{ref}}.
\end{aligned}
\end{equation}
Herein, we have
\begin{equation}\label{eq:diff+rhs}
{\bs A}({\hat x},{\bs y})\isdef({\bs J}^\intercal{\bs J})^{-1}({\hat x},{\bs y})\det{\bs J}({\hat x},{\bs y}),\quad\hat{f}({\hat x},{\bs y})\isdef (f\circ{\bs V})({\hat x},{\bs y})
\det{\bs J}({\hat x},{\bs y})
\end{equation}
and 
\[ \hat{u}({\hat x},{\bs y})\isdef (u\circ{\bs V})({\hat x},{\bs y}).
\]

\begin{remark}
The uniformity condition in \eqref{eq:unifVfield} implies that the functional determinant 
\(\det {\bs J}({\hat x},{\bs y})\) in \eqref{eq:diff+rhs}
is either uniformly positive or negative, see \cite{HPS16} for the details. 
We shall assume without loss of generality \(\det {\bs J}({\hat x},{\bs y})>0\)
and hence \(|\det {\bs J}({\hat x},{\bs y})|=\det {\bs J}({\hat x},{\bs y})\), i.e.\ we may just drop the modulus.
More precisely, due to \eqref{eq:unifVfield}, we can bound the determinant according to
\[
0<\frac{1}{C_{\on{uni}}^d}\leq\det {\bs J}({\hat x},{\bs y})\leq C_{\on{uni}}^d<\infty
\]
for every \({\hat x}\in D_{\on{ref}}\) and almost every \({\bs y}\in\Gamma_\infty\). 
In addition, all singular values of  
 \({\bs J}^{-1}({\hat x},{\bs y})\) are bounded from below by \(C_{\on{uni}}^{-1}\) and from above by \(C_{\on{uni}}\).
From this, we obtain the bound
\begin{equation}\label{eq:uniEllipticity}
0<\check c:=\frac{1}{C_{\on{uni}}^{d+2}}\leq\| {\bs A}({\hat x},{\bs y})\|_2\leq C_{\on{uni}}^{d+2}=:\hat c<\infty
\end{equation}
for every \({\hat x}\in D_{\on{ref}}\) and almost every \({\bs y}\in\Gamma_\infty\). Hence, the transported model problem is uniformly
elliptic.
\end{remark}

We conclude this section by summarizing the regularity results for \({\bs A},\hat{f},\hat{u}\), cp.~\eqref{eq:parmodProb}, with respect
to the parameter \({\bs y}\in\Gamma_\infty\) from \cite{HPS16}.
For this, denote by $\mcF$ the set of finitely supported multi-indices
\begin{equation*}
\mcF \isdef \mset{\mu\in\bbN_0^\infty}{\abs{\supp\mu} <\infty}\quad\text{where}\quad \supp\mu \isdef \mset{m\in\bbN}{\mu_m\neq 0}.
\end{equation*}

\begin{theorem}
\label{thm:regularity}
Let the right-hand side \(f\) from \eqref{eq:modProb} satisfy \[\|\partial^{\bs\alpha}_{x}f\|_{L^\infty(\mathcal{D})}\leq c_f{\bs\alpha}!\rho^{-|\bs\alpha|}\]
for some constants \(c_f,\rho>0\) and $\alpha\in\mcF$.
Then, for every \({\bs\alpha}\in\mathcal F\) it holds
\begin{align*}
\big\|\partial^{\bs\alpha}_{\bs y}{\bs A}({\bs y})\big\|_{L^\infty(D_{\on{ref}};\mathbb{R}^{d\times d})}
&\leq C|\balpha|!c^{|\bs\alpha|}{
\bs\gamma
}^{\bs\alpha},\\
\big\|\partial^{\bs\alpha}_{\bs y}\hat{f}({\bs y})\big\|_{L^\infty(D_{\on{ref}})}
&\leq C|\balpha|!c^{|\bs\alpha|}{
\bs\gamma
}^{\bs\alpha},\\
\big\|\partial^{\bs\alpha}_{\bs y}\hat{u}({\bs y})\big\|_{H^1(D_{\on{ref}})}
&\leq C|\balpha|!c^{|\bs\alpha|}{
\bs\gamma
}^{\bs\alpha},\\
\end{align*}
for some constants \(C>0\), which depend on \(c_f,\rho,C_{\on{uni}},d,D_{\on{ref}},\|\bs\gamma\|_{\ell^1}\)
but are independent of the multi-index \(\bs\alpha\).
\end{theorem}

%% file: galerkin.tex
\section{Adaptive Galerkin discretisation}
\label{sec:adaptive galerkin}

In this section we describe the Galerkin discretization of the considered random PDE~\eqref{eq:parmodProb} in a finite dimensional subspace $\mathcal V_N\subset\mathcal V=\mathcal X \otimes \mathcal Y$.
Determined by the elliptic problem type with homogeneous boundary condition, we assume $\mcX = H^1_0(D_{\on{ref}})$ is discretized by a first order Lagrange FE basis on a mesh representing $D_{\on{ref}}$. 
Moreover, the randomness is modelled in a truncated version of $\mcY =  L^2(\Gamma_\infty, \pi_\infty)$ and represented by Legendre chaos polynomials orthonormal with respect to the joint probability measure $\pi_\infty$ associated with the parameter $\bs y$.
Consequently, $\mcV_N$ with norm $\norm{v}^2_\mcV=\bbE_{\pi_M}[\norm{v}^2_\mcX]$ is spanned by the respective tensor basis.
Moreover, the residual based a posteriori error estimator of \cite{EGSZ1,EGSZ2} is recalled for the problem at hand.

For efficient computations of the Galerkin projection and the error estimator, the resulting system with inhomogeneous coefficient and right-hand side~\eqref{eq:diff+rhs} is represented in the tensor train format as presented in Section~\ref{sec:tensor trains}.


\subsection{Parametric and deterministic discretization}
\label{sec:discretization}

To determine a multivariate polynomial basis of $\mathcal Y$, we first define the full tensor index set of order $M\in\bbN$ and maximal degree $(d_1,\ldots, d_M)\in\bbN^M$ by
\begin{equation*}
\Lambda \isdef \{(\mu_1,\ldots,\mu_M,0,\ldots) \in \mathcal F : \mu_m = 0,\ldots, d_m-1, \ m = 1,\ldots,M \}
\end{equation*}

For any such subset $\Lambda\subset\mcF$, we define $\supp\Lambda\isdef\bigcup_{\mu\in\Lambda} \supp\mu \subset \bbN$.
Let $(P_n)_{n=0}^\infty$ denote a basis of $L^2([-1,1])$, orthogonal with respect to the Lebesgue measure, consisting of Legendre polynomials $P_n$ of degree $n\in\bbN_0\
$ on $\bbR$.
Moreover, to obtain a finite dimensional setting, we define the truncated parameter domain $\Gamma_M \isdef [-1, 1]^M$ and probability measure $\pi_M \isdef \bigotimes_{m=1}^M \pi_m$.
By tensorization of the univariate polynomials, an orthogonal basis of $L^2(\Gamma_M)=\bigotimes_{m=1}^M L^2([-1, 1])$ is obtained.
Then, for any multi-index $\mu\in\mcF$, the tensor product polynomial $P_\mu \isdef \bigotimes_{m=1}^\infty P_{\mu_m}$ in $\bs y\in\Gamma_M$ is expressed as the finite product
\begin{equation*}
P_\mu(\bs y) = \prod_{m=1}^ \infty P_{\mu_m}(y_m) = \prod_{m\in\supp \mu}P_{\mu_m}(y_m).
\end{equation*}
Assuming that $\pi_M=\mathcal{U}(\Gamma_M)$, after suitable rescaling we can consider $(P_\mu)_{\mu\in\mcF}$ as an orthonormal basis of $L^2(\Gamma_M, \pi_M) = \bigotimes_{m=1}^M L^2([-1, 1], \frac{1}{2}\lambda)$, where $\lambda$ denotes the Lebesgue measure and hence $\frac{1}{2}\lambda$ is the uniform measure on $[-1, 1]$, see~\cite{SG}.

A discrete subspace of $\mcX$ is given by the conforming finite element space
$\mcX_p(\mcT)\isdef \mathrm{span} \{\varphi_i\}_{i=1}^N \subset \mcX$ of degree $p\geq 0$ on some simplicial regular mesh $\mcT$ of the domain $D_{\on{ref}}$ with the set of faces $\mcS$ (i.e.\ edges for $d=2$) and basis functions $\varphi_i$.
For a convenient presentation, we denote the piecewise constant basis functions of $\mcX_0(\mcT)$ by $\{\psi_\ell\}_{\ell=1}^{N_0}$, where $N_0=\dim \mcX_0$ is the number of elements in $\mcT$.
In order to circumvent complications due to an inexact approximation of boundary values, we assume that $D_{\on{ref}}$ is a polytope.
By denoting $P_p(\mcT)$ the space of piecewise polynomials of degree $p\geq 0$ on $\mcT$, the assumed FE discretization with Lagrange elements then satisfies $\mcX_p(\mcT)\subset P_p(\mcT)\cap C(\overline{\mcT})$.
For any element $T\in\mcT$ and face $S\in\mcS$, we set the entity sizes $h_T\isdef \diam T$ and $h_S \isdef \diam S$.
Let $n_S$ denote the exterior unit normal on any face $S$.
The jump of some $\chi\in H^1(D_{\on{ref}};\bbR^d)$ on $S = \overline{T_1}\cap\overline{T_2}$ in normal direction is then defined by
\begin{equation}
\jump{\chi}_S \isdef \chi\vert_{T_1}\cdot n_S - \chi\vert_{T_2}\cdot n_S.
\end{equation}
By $\omega_T$ and $\omega_S$ we denote the element and facet patches defined by the union of all elements which share at least a vertex with $T$ or $S$, respectively.
Consequently, the Cl\'ement interpolation operator $I\colon \mcX \to \mcX_p(\mcT)$ satisfies, respectively for $T\in\mcT$ and $S\in\mcS$,
\begin{align}
\label{eq:clement}
    \norm{(\id-I) v}_{L^2(T)} \leq c_{\mathcal T} h_T \abs{v}_{\mcX,\omega_T},\qquad
    \norm{(\id-I) v}_{L^2(S)} &\leq c_{\mathcal S} h_S^{1/2} \abs{v}_{\mcX,\omega_S},
\end{align}
where the seminorms $\abs{\;\cdot\;}_{\mcX,\omega_T}$ and $\abs{\;\cdot\;}_{\mcX,\omega_S}$ are the restrictions of $\norm{\;\cdot\;}_{\mcX}$ to $\omega_T$ and $\omega_S$, 

The fully discrete approximation space subject to some mesh $\mathcal T$ with FE order $p\geq 0$ and active set $\Lambda$ with $\abs{\Lambda}<\infty$ is given by
\begin{equation}
\mcV_N\isdef \mcV_N(\Lambda,\mcT,p)\isdef \bigg\{v_N(x,y) = \sum_{\mu\in\Lambda} v_{N,\mu}(x)P_\mu(y): \ v_{N,\mu}\in \mcX_p(\mcT)\bigg\}, 
\end{equation}
and it holds $\mcV_N\subset\mcV$.
We define a tensor product interpolation operator $\mcI\colon L^2(\Gamma_\infty,\pi_\infty;\mcX) \to \mcV_N$
for $v = \sum_{\mu \in \mcF} v_\mu P_\mu \in \mcV = L^2(\Gamma_\infty,\pi_\infty;\mcX)$ by setting
\begin{equation}
\mcI v \isdef \sum_{\mu \in \Lambda} (I v_\mu) P_\mu.
\end{equation}
For $v \in \mcV$ and all $T \in \mcT,\ S \in \mcS$, it holds
\begin{align}
\label{eq:tensor-clement}
\norm{(\id - \mcI)v}_{L^2(\Gamma_\infty,\pi_\infty;L^2(T))} &\leq c_{\mathcal T} h_T \abs{v}_{\mcV,\omega_T}, \\
\norm{(\id- \mcI) v}_{L^2(\Gamma_\infty,\pi_\infty;L^2(S))} &\leq c_{\mathcal S} h_S^{1/2} \abs{v}_{\mcV,\omega_S},
\intertext{where}
\notag
\abs{v}_{\mcV,\omega_T}^2 \isdef \int_{\Gamma_\infty} \abs{v(y)}_{\mcX,\omega_T}^2 \dx{\pi_\infty}(y),&\qquad  \abs{v}_{\mcV,\omega_S}^2 \isdef \int_{\Gamma_\infty} \abs{v(y)}_{\mcX,\omega_S}^2 \dx{\pi_\infty}(y).
\end{align}

\subsection{Random field discretisation}
\label{sec:field discretization}
In this paragraph, we highlight the special structure of the random field discretization.
We aim at an efficient way to discretize the transformed and parametrized random fields~\eqref{eq:diff+rhs} in terms of the piecewise constant finite element functions $\{\psi_i\}_{i=1}^{N_0}\subset\mcX_0\subset\mcX$, pointwise for every $\bs y\in\Gamma_M$.

In \cite{MM18}, it has been shown how the random vector field \eqref{eq:ParVfield} can efficiently be approximated by means of finite elements. 
This results in a truncated representation with \(M\in\mathbb{N}\) terms of the form
\[
{\bs V}_h(\hat{x},{\bs y})=
\hat{x}+\sum_{m=1}^My_m{\bs V}_{m,h}(\hat{x})=\hat{x}+\sum_{m=1}^My_m\sum_{i=1}^{d}\sum_{k=1}^N c_{i,k,m}\varphi_k(\hat{x}){\bs e}_i,
\]
where \({\bs e}_1,\ldots,{\bs e}_d\) denotes the canonical basis of \(\mathbb{R}^d\), \(\varphi_1,\ldots,\varphi_N\) is a basis for \(\mcX_1(\mcT)\) and \(c_{i,k,m}\in\mathbb{R}\) are the coefficients in the basis representation of \({\bs V}_{m,h}\).
The length \(M\) of this expansion depends on the desired approximation error of the random field, which can be rigorously controlled in terms of operator traces, see \cite{HPS14,ST2}.

For the corresponding Jacobian, we obtain
\[
{\bs J}_h(\hat{x},{\bs y})={\bs I}+\sum_{m=1}^My_m{\bs V}_{m,h}'(\hat{x})={\bs I}+\sum_{m=1}^My_m\sum_{i=1}^{d}\sum_{k=1}^N c_{i,k,m}{\bs e}_i\big(\nabla_{\hat{x}}\varphi_k(\hat{x})\big)^\intercal.
\]
More explicitly, the Jacobians \({\bs V}_{m,h}'(\hat{x})\) are given by
\[
{\bs V}_{m,h}'(\hat{x})=\sum_{k=1}^N\begin{bmatrix}c_{1,k,m}\partial_1\varphi_k(\hat{x})&\cdots &c_{1,k,m}\partial_d\varphi_k(\hat{x})\\
\vdots & \ddots & \vdots\\
c_{d,k,m}\partial_1\varphi_k(\hat{x})&\cdots &c_{d,k,m}\partial_d\varphi_k(\hat{x})
\end{bmatrix}.
\]
Since \(\partial_i\varphi_k(\hat{x})\), \(i=1,\ldots,d\), \(k=1,\ldots,N\) are piecewise constant functions, we can represent \({\bs V}_{m,h}'\) in an element based fashion according to
\[
{\bs V}_{m,h}'=\sum_{\ell=1}^{N_0}\begin{bmatrix}\tilde{c}_{\ell,m,1,1} & \cdots & \tilde{c}_{\ell,m,1,d}\\
\vdots & \ddots & \vdots \\
\tilde{c}_{\ell,m,d,1} & \cdots & \tilde{c}_{\ell,m,d,d}
\end{bmatrix}\psi_\ell(\hat{x})\defis\sum_{\ell=1}^{N_0}{\bs C}_{\ell,m}\psi_\ell(\hat{x}),
\]
where \(\psi_\ell\) denotes the characteristic function on the element \(T_\ell\in\mathcal{T}\) and \(\tilde{c}_{\ell,m,i,j}\in\mathbb{R}\) 
are the corresponding coefficients. 
Hence, we end up with a piecewise constant representation of \({\bs V}_h'\), which reads
\[
{\bs J}_h(\hat{x},{\bs y})={\bs I}+\sum_{\ell=1}^{N_0}\bigg(\sum_{m=1}^M{\bs C}_{\ell,m}y_m\bigg)\psi_\ell(\hat{x}).
\]
From this representation, it is straightforward to calculate \(\det{\bs J}_h({\hat{x}},{\bs y})\) for a given \({\bs y}\in\Gamma_M\),
also in an element based fashion. Having \({\bs V}_h({\hat{x}},{\bs y})$, ${\bs J}_h({\hat{x}},{\bs y})$, $\det{\bs J}_h({\hat{x}},{\bs y})\)
at our disposal, it is then easy to evaluate \({\bs A}({\hat x},{\bs y})\) and \(\hat{f}({\hat x},{\bs y})\), as well.

This procedure can be extended to the general case of order $p>0$ ansatz functions for the random vector field~\eqref{eq:ParVfield}, resulting in an order $p-1$ approximation of the desired quantities in~\eqref{eq:diff+rhs}.

\subsection{Variational formulation}
\label{sec:variational formulation}

Using the transformation in~\eqref{eq:diff+rhs}, the weak formulation of the model problem~\eqref{eq:parmodProb} reads: find $u\in\mcV$, such that for all $v\in\mcV$ there holds
\begin{multline}
  \label{eq:variational formulation}
  \int_{\Gamma_\infty} \int_{D_{\on{ref}}}\bs A(\hat x, \bs y) \nabla_{\hat x}u(\hat x, \bs y) \cdot\nabla_{\hat x}v(\hat x, \bs y) \dx{\hat x} \dx{\pi_\infty}(\bs y) =\\
  	\langle \mcA(u), v\rangle=
 \int_{\Gamma_\infty} \int_{D_{\on{ref}}} \hat{f}(\hat x, \bs y) v(\hat x, \bs y) \dx{\hat x} \dx{\pi}(\bs y).
\end{multline}
This characterizes the operator $\mcA:L^2(\Gamma_\infty,\pi_\infty;V)\to L^2(\Gamma_\infty,\pi_\infty;V^\ast)$, which gives rise to the energy norm $\norm{v}^2_\mcA\isdef\langle \mcA(v),v\rangle$.
Employing the finite dimensional spaces of the previous section leads to the discrete weak problem: find $u=\sum_{\mu\in\Lambda}\sum_{i=1}^N U(i, \mu)\varphi_iP_\mu\in\mcV_N$, such that for all $i'=1, \ldots, N$ and $\alpha'\in\Lambda$
\begin{equation*}
\sum_{\alpha\in\Lambda}\sum_{i=1}^N\mathbf{L}(i, \alpha, i', \alpha')U(i', \alpha') = \mathbf{F}(i', \alpha').
\end{equation*}
Here, we define the discrete linear operator
\begin{equation}
  \label{eq:discrete-operator}
  \mathbf{L}(i, \alpha, i', \alpha') \isdef \int_{\Gamma_M}\int_{D_{\on{ref}}}\bs A(\hat x, \bs y) \nabla_{\hat x}\varphi_i(\hat{x})P_\alpha(\bs y)\nabla_{\hat x}\varphi_{i'}(\hat{x})P_{\alpha'}(\bs{y}) \dx{\hat{x}}\dx{\pi_M}(\bs y)
\end{equation}
and the discrete right-hand side
\begin{equation*}
  \mathbf{F}(i', \alpha') \isdef \int_{\Gamma_M}\int_{D_{\on{ref}}}\hat f(\hat x,\bs y) \varphi_{i'}(\hat{x}) P_{\alpha'}(\bs y) \dx{\hat x}\dx{\pi_M}.
\end{equation*}


\subsection{Residual based a posteriori error estimates}
\label{sec:error estimator}

In the following, we recall the residual based a posteriori error estimator derived in~\cite{EGSZ1,EGSZ2}, adopted for the problem at hand.
An efficient reformulation in the tensor train format is postponed to Section~\ref{sec:tensor trains}.
The basis for the estimator is the residual $\mcR(w_N) \in L^2(\Gamma_\infty, \pi_\infty;\mcX^\ast) = \mcV^\ast$ with respect to some $w_N\in\mcV_N$ and the solution $u\in\mcV$ of~\eqref{eq:parmodProb} given by
\[
\mcR(w_N) \isdef \mcA(u-w_N) = \hat f - \mcA(w_N).
\]
It has an $L^2(\Gamma_\infty, \pi_\infty)$-convergent expansion in $(P_\nu)_{\nu\in\mcF}$ given by
\[
\mcR(w_N) = \sum_{\nu\in\mcF} r_\nu(w_N)P_\nu,
\]
with coefficients $r_\nu\in\mcX^\ast$ characterized by
\begin{equation}
\label{eq:rnu}
\langle r_\nu, v \rangle = \int_{D_{\on{ref}}} \hat f_\nu v - \sum_{(\mu,\kappa)\in\Upsilon_\nu} \bs A_\mu\nabla_{\hat{x}} w_{N,\kappa}\beta_{\mu, \kappa}^{\nu} \cdot \nabla_{\hat{x}} v\dx{\hat x}\qquad \forall v\in\mcX.
\end{equation}
Here, $\hat f_\nu$, $\bs A_\mu$ and $w_{N, \kappa}$ denote the coefficients in the Legendre chaos expansion of $\hat f=\sum_{\nu\in\mcF}\hat f_\nu P_\mu$, $\bs A=\sum_{\mu\in\mcF}\bs A_\mu P_\mu$ and $w_N = \sum_{\kappa\in\mcF}w_{N, \kappa} P_\kappa$ and 
$$\Upsilon_\nu\isdef\bigg\{(\mu,\kappa)\in\mcF\times\Lambda : \beta_{\mu, \kappa}^{\nu} \isdef \int_{\Gamma_\infty} P_\nu(\bs y) P_\mu(\bs y) P_\kappa(\bs y) \dx{\pi_\infty}(\bs y)\neq 0\bigg\}$$ 
is the $\nu$-relevant triple product tuple set.

We recall a central theorem from~\cite{EGSZ1}, which enables the derivation of an error bound based on an approximation $w_N$ of the Galerkin projection $u_N$ of the solution $u$ in the energy norm.
\begin{theorem}
\label{thm:error bound}
Let $\mcV_N\subset\mcV$ be a closed subspace and $w_N\in\mcV_N$, and let $u_N\in\mcV_N$ denote the $\mcA$ Galerkin projection of $u\in\mcV$ onto $\mcV_N$. Then, for some $c_\mcA,c_\mcI>0$, it holds
\begin{align*}
\norm{u-w_N}_\mcA^2 \leq \check c^2 \left( \sup_{v\in\mcV} \frac{\abs{\langle \mcR(w_N),(\id-\mcI)v\rangle}}{\norm{v}_\mcV} + c_\mcI\norm{u_N-w_N}_\mcA \right)^2 + \norm{u_N-w_N}^2_\mcA.
\end{align*}
\end{theorem}
\begin{remark}
The constant $c_\mcI$ is related to the Cl\'ement interpolation operator in $V$ and $\check c$ stems from the spectral equivalence such that $\norm{v}_\mcA\geq \check c\norm{v}_\mcV$.
We refer to~\cite{EGSZ1} for further details.
\end{remark}
\begin{remark}
\label{rem:finite poly dim}
We henceforth assume that the data $\hat f$ and $\bs A$ are exactly expanded in a finite set $\Delta$ with $\Lambda\subset\Delta\subset\mcF$, i.e. with the approximation, there is no significant contribution from the neglected modes $\mcF\setminus\Delta$.
The residual can then be split into approximation and truncation contributions
\begin{equation}
\label{eq:residual split}
\mcR(w_N) = \mcR_\Lambda(w_N) + \mcR_{\Delta\setminus\Lambda}(w_N),
\end{equation}
where $\mcR_\Xi$ denotes the restriction of the expansion to the set $\Xi\subset\mcF$.
Computable upper bounds for the two residual terms and the algebraic error $\norm{u_N-w_N}_\mcA$ are recalled in the following.
\end{remark}

For any discrete $w_N\in\mcV_N$, we define the following error estimators in analogy to the presentation in~\cite{EGSZ1,EGSZ2} and~\cite{eigel2018adaptive}:

\begin{itemize}
\item A deterministic residual estimator for $\mcR_\Lambda$ steering the adaptivity of the mesh $\mcT$ is given by
\begin{equation}
  \label{eq:det estimator intro}
  \eta(w_N)^2 \isdef \sum_{T\in\mathcal{T}} \eta_T(w_N)^2 + \sum_{S\in\mathcal{S}} \eta_S(w_N)^2,
\end{equation}
with volume contribution for any $T\in\mcT$ 
\begin{equation}
\label{eq:est-cont}
\eta_T(w_N) \isdef \bigg\|\sum_{\nu\in\Lambda} \left(\hat f_\nu - \div\left(\sum_{(\mu, \kappa)\in\Upsilon_\nu} \bs A_\mu\nabla_{\hat{x}} w_{N, \kappa}\beta_{\mu, \kappa}^\nu \right)\right)P_\mu
\bigg\|_{L^2(\Gamma_M, \pi_M; L^2(T))}
\end{equation}
and facet contribution for any $S\in\mcS$
\begin{equation}
\eta_S(w_N) \isdef \bigg\|\sum_{\nu\in\Lambda} \bigjump{\sum_{(\mu, \kappa)\in\Upsilon_\nu} \bs A_\mu\nabla_{\hat{x}} w_{N, \kappa}\beta_{\mu, \kappa}^\nu}_S P_\nu\bigg\|_{L^2(\Gamma_M, \pi_M; L^2(S))}.
\end{equation}

\item The stochastic truncation error estimator stems from splitting the residual in~\eqref{eq:residual split}, while considering the inactive part over $\Delta\setminus\Lambda$.
It is possible to construct the estimator, as in the deterministic case, for every element of the triangulation and consider different mesh discretisations for every stochastic multi-index.
Since we want to focus on a closed formulation and avoid technical details, the stochastic estimator is formulated on the whole domain $D_{\on{ref}}$.
Nevertheless, for more insight we introduce a collection of suitable tensor sets, which indicate the error portion of every active stochastic dimension $m=1,\ldots,M$ (in fact, we could even consider $m>M$),
\begin{align}
  \label{eq:Delta_n}
  \Delta_m \isdef \{\mu\in\mathcal{F}\;\vert\;& \mu_j=0, \ldots, d_j-1,\; j=1, \ldots, M, \notag \\ & j\neq m, \mu_m = d_m, \mu_k=0, k > M\}.
\end{align}
Then, for every $w_N\in\mathcal{V}_N$, the stochastic tail estimator on $\Delta_m$ is given by
\begin{equation}
  \label{eq:sto estimator intro}
  \zeta_m(w_N) \isdef \bigg\|\sum_{\nu\in\Delta_m} \zeta_\nu(w_N) P_\nu\bigg\|_{L^2(\Gamma_M, \pi_M; L^2(D_{\on{ref}}))}, 
\end{equation}
where we define for every multi index $\nu\in\mcF$ the residual portion
\begin{equation}
\zeta_\nu \isdef \hat f_\nu - \div\left(\sum_{(\mu, \kappa)\in\Upsilon_\nu} \bs A_\mu\nabla_{\hat{x}} w_{N, \kappa}\beta_{\mu, \kappa}^\nu\right)\in\mcX^\ast.
\end{equation}
The collection of sets $\{\Delta_n\}_{n=1}^M$ is beneficial in the adaptive refinement procedure but it does not cover the whole stochastic contributions of the residual.
For this, we need to compute the global stochastic tail estimator over $\mcF\setminus \Lambda$
\begin{equation}
  \label{eq:global sto estimator}
  \zeta(w_N) \isdef \norm{\sum_{\nu\in\mcF\setminus \Lambda} \zeta_\nu(w_N) P_\nu}_{L^2(\Gamma_M, \pi_M; L^2(D_{\on{ref}}))}
\end{equation}
 which incorporates an infinite sum that becomes finite due to remark~\ref{rem:finite poly dim}.

\item The algebraic error denotes the distance of $w_N$ to the $\mcV_N$ best approximation $u_N$.
In particular, this distance can e.g. occur due to an early termination of an iterative solver or an restriction to another solution manifold $\mathcal{M}\subset \mcV_N$.
This error can be bounded by
\begin{equation}
  \norm{u_N - w_N} \leq \iota(w_N),
\end{equation}
where
\[
\iota(w_N) \isdef \norm{(\bs L W - \bs F)\bs H^{-1/2}}_F.
\]
Here, $W\in\bbR^{N, d_1,\ldots, d_M}$ denotes the coefficient tensor of $w_N\in\mathcal{V}_N$, $\bs L$ is the discrete operator from~\eqref{eq:discrete-operator} and $\norm{\cdot}_F$ is the Frobenius norm.
Note that the rank-1 operator $\bs H$ is a base change operator to orthonormalize the physical basis functions, i.e.\
\begin{equation}
  \bs H \isdef H_0 \otimes I \otimes \dots \otimes I, \quad H_0(i, i') = \int_{D_{\on{ref}}} \nabla_{\hat{x}}\varphi_i(\hat{x}) \nabla_{\hat{x}}\varphi_{i'}(\hat{x})\mathrm{d}\hat{x}.
\end{equation}
\end{itemize}

The combination of these estimators in the context of Theorem~\ref{thm:error bound} yields an overall bound $\Theta$ for the energy error similar to the references~\cite{EGSZ1, EGSZ2, EPS} and~\cite{eigel2018adaptive}
\begin{corollary}
\label{cor:error estimator}
For any $w_N\in\mathcal{V}_N$, the solution $u\in\mathcal{V}$ of the model problem~\eqref{eq:modProb} and the Galerkin approximation $u_N\in\mathcal{V}_N$ in~\eqref{eq:variational formulation}, there exists constants $c_\eta, c_\zeta, c_\iota>0$ such that it holds
\begin{align}
  \label{eq:error estimator}
  \norm{w_N - u}_\mcA^2 \leq \Theta^2 \isdef \bigl( c_\eta \eta(w_N) + c_\zeta \zeta(w_N) + c_\iota \iota(w_N)\bigr)^2 + \iota(w_N)^2.
\end{align}
\end{corollary}

\begin{remark}
Observing the residual in~\eqref{eq:rnu} it becomes clear that the derived error estimators suffer from the ``curse of dimensionality'' and are hence not computable for larger problems. However, the hierarchical low-rank tensor representation introduced in the next section alleviates this substantial obstacle and makes possible the adaptive algorithms described in Section~\ref{sec:adaptive algorithm}.
\end{remark}

%% file: tensorTrains.tex
\section{Tensor trains}
\label{sec:tensor trains}

The inherent tensor structure of the involved Bochner function space \;$\mathcal{V} = \mathcal{X}\bigotimes_{m=1}^M L^2([-1, 1], \pi_m)$ and the corresponding finite dimensional analogue $\mathcal{V}_N$ motivates the use of hierarchical tensor formats which aim at an implicit model order reduction, effectively breaking the curse of dimensionality in case of low-rank approximability of operator and solution.

A representative $v\in\mathcal{V}_N$ can be written as 
\begin{equation}
  \label{eq:discrete function}
  v(\hat x, \bs{y}) = \sum_{i=1}^N\sum_{\mu\in\Lambda} V(i, \mu)\varphi_i(\hat{x})P_{\mu}(\bs{y}),
\end{equation}
where $V\in \mathbb{R}^{N, d_1,\ldots, d_M}$ is a high dimensional tensor containing for example the projection coefficients 
$$V(i, \mu) = \mathbb{E}_{\pi_M}\bigg[\int_{D_{\on{ref}}} v(\hat{x}, \cdot) \varphi_i(\hat{x})\mathrm{d}x P_\mu(\cdot)\bigg].$$
Setting $d=\max\{d_1,\ldots, d_M\}$, the storage cost of $V$ is $\mathcal{O}(N d^M)$, which grows exponentially with the number of dimensions $M\in\mathbb{N}$ in the stochastic parameter space.
To alleviate this major problem for numerical methods, we impose a \emph{low-rank} assumption on the involved objects and introduce a popular tensor format as follows. 

A tensor $V\in\mathbb{R}^{N, d_1,\ldots, d_M}$ is called in \emph{tensor train} (TT) format if every entry can be represented as the result of a matrix-vector multiplication of the form
\begin{equation}
  \label{eq:TT format}
  V(i, \mu_1,\ldots, \mu_M) = \sum_{k_0=1}^{r_0}\dots \sum_{k_{M-1}=1}^{r_{M-1}}V_0(i, k_0) \prod_{m=1}^M V_m(k_{m-1}, \mu_m, k_m).
\end{equation}
To simplify notation, set $r_M=1$. 
If the vector $\mathbf{r}=(r_0, \ldots, r_M)$ is minimal in some sense, we call $\mathbf{r}$ the \emph{TT-rank} and~\eqref{eq:TT format} is the \emph{TT}-decomposition of $V$.
It can be observed that the complexity of $V$ now depends only linearly on the number of dimensions, namely $V=\mathcal{O}(dM\max\{\mathbf{r}\}^2)$.
In~\cite{oseledets2011tensor, kolda2009tensor} it was shown that many functions in numerical applications admit a low-rank representation. 

Given the full tensor description of $V$, one could compute the tensor train representation by a \emph{hierarchical singular value decomposition} (HSVD) as described in~\cite{grasedyck2010hierarchical}.
However, this is usually unfeasible due to the high dimensionality of $V$ or because it is known only implicitly.
In that case, the utilization of high dimensional interpolation or regression algorithm is advisible, see e.g.~\cite{oseledets2010tt, grasedyck2015variants}.

In this work, we rely on a TT reconstruction approach and employ it to obtain the representation of the transformed coefficient function and the right-hand side~\eqref{eq:diff+rhs}.
Opposite to an explicit (intrusive) discretisation of the linear system in tensor format as e.g. carried out in~\cite{eigel2018adaptive,MR2783199}, the reconstruction method relies on a set of random samples of the solution.
The non-intrusive algorithm used in the numerical experiments is described in~\cite{eigel2019non}.
Similar ideas were presented in~\cite{oseledets2010tt,dolgov2015polynomial,dolgov2017hybrid}, where a tensor cross approximation was used for the construction of the algebraic system.
In contrast to the tensor reconstruction, a selective sampling of strides in the tensor has to be available to perform a cross approximation.
Consider~\cite{grasedyck2013literature} for a survey on the topic of low-rank approximations methods.

To sketch the reconstruction approach, we assume a set $\{y^{(k)}\}_{k=1}^K$ of $K$ parameter realisations and corresponding measurements of a function $\{v(\cdot, y^{k})\in\mcX\}_{k=1}^K$
\[
{\bf b}=(v(\cdot,y^{(1)}),\ldots,v(\cdot,y^{(K)}))^\intercal.
\]
Recall that $N$ is the dimension of the finite element space.
We define a linear measurement operator $\hat\mcA:\bbR^{N\times d_1\times\cdots\times d_M}\to\bbR^{NK}$ acting on a tensor $W\in\bbR^{N\times d_1\times\cdots\times d_M}$ by
\begin{equation*}
\hat\mcA(W)_{j,k} = (W\circ_M({\bf\xi}_1^{(k)}\otimes\cdots\otimes{\bf\xi}_M^{(k)}))[j],\quad 1\leq k\leq K, 1\leq j\leq N,
\end{equation*}
with a contraction $\circ_M$ over the $M$ stochastic modes and
\begin{equation*}
{\bf\xi}_m^{(k)} := \left(P_1(y_m^{(k)}) \cdots P_{d_m}(y_m^{(k)}) \right)^\intercal.
\end{equation*}
The reconstruction problem is to find a tensor $W$ with minimal TT-rank such that $\hat\mcA(W)={\bf b}$.
Details in particular of the numerical solution algorithm of the optimisation problem by an Alternating Steepest Descent (ASD) can be found in~\cite{eigel2019non}.

\subsection{Galerkin discretization in tensor train format}

In the following, we assume an acessible tensor representations of the right-hand side $\hat f$ and the coefficient function $\bs A$.
To make this more precise, we denote the low-rank approximations of e.g.\ $\hat f$ in~\eqref{eq:diff+rhs} by
\begin{equation}
  \label{eq:tt-rhs}
	\hat f(\hat x, \bs y) \approx f^{\mathrm{TT}}(\hat{x}, \bs{y}) = \sum_{\mu\in\Lambda_f}\sum_{i=1}^{N_0} F(i, \mu)\psi_i(\hat{x})P_{\mu}(\bs{y}),
\end{equation}
where $F$ admits a TT representation of rank $\mathbf{r}^f$ and $\Lambda_f$ is a tensor multi-index set with local dimension cap $\mathbf{d}^f=(d_1^f, \ldots, d_M^f)$.
Analogously, every component of the symmetric matrix coefficient 
\begin{equation}
  \label{eq:coef-mat}
  \bs{A}(\hat{x}, \bs{y}) = \begin{bmatrix}
    \bs{a}_{1, 1}(\hat{x}, \bs{y}) &  \cdots &   \bs{a}_{1, d}(\hat{x}, \bs{y}) \\ 
    \vdots & \ddots & \vdots \\
        \bs{a}_{d, 1}(\hat{x}, \bs{y}) &  \cdots  & \bs{a}_{d, d}(\hat{x}, \bs{y})
  \end{bmatrix}
\end{equation}
is approximated by $a^{\mathrm{TT}}_{i, j}$, $i,j\in\{1, \ldots, d\}$ as in~\eqref{eq:tt-rhs} with TT-ranks $\mathbf{r}^{i,j}$.
Here, the order three component tensors in the TT-representation of the approximated matrix entry $\bs A^{\mathrm{TT}}\vert_{i, j} = a_{i,j}^{\mathrm{TT}}$ are denoted by $\{a_{i,j}^{m}\}_{m=0}^M$.

\begin{remark}
\label{rem:ttdim}
Since, for the coefficient, the TT reconstruction is carried out for every matrix entry in~\eqref{eq:coef-mat}, the local dimensions $d^{i, j}=(d_1^{i, j}, \ldots, d_M^{i, j})$ and tensor ranks can vary among those $d^2$ tensor trains.
Here, we assume that every approximation has the same local dimensions and the tensor multi-index set covering those indices is denoted by $\Xi\subset\mcF$, possibly different from (but larger than) the solution active set $\Lambda$.
As stated in~\cite{eigel2018adaptive}, it is beneficial (and in fact necessary) to chose $\Xi$ such that for all $\mu\in\Lambda$ also $2\mu=(2\mu_1,\ldots, 2\mu_M, 0,\ldots)\in\Xi$.
Due to the orthogonality of the polynomial basis $\{P_\nu\}$, this feature guarantees a well-posed discrete problem since additional approximations are avoided and enables quasi-optimal convergence rates of the Galerkin method.
\end{remark}

%
On $\mathcal{V}_N$, the Galerkin operator resulting from the transformed weak problem in TT format is given as the sum of $d^2$ TT operators such that for all $i, i'=1, \ldots, N$, and $\alpha,\alpha'\in\Lambda$,
\begin{equation}
  \label{eq:TTop sum}
  \mathbf{L}(i, \alpha, i', \alpha') = \sum_{j=1}^{d^2}\mathbf{L}_j(i, \alpha, i', \alpha'),
\end{equation}
each corresponding to one addend of the resulting matrix-vector product in~\eqref{eq:discrete-operator}. 

In the following, we illustrate the explicit construction of the TT operator for the term $\mathbf{L}_1$.
By denoting $\nabla^{i}g$ the $i$-th component of the gradient of a function $g$, for the first low-rank approximated bilinear form addend one obtains
\begin{equation}
  \mathbf{L}_1(i, \alpha, i', \alpha') \approx \int_{\Gamma_M}\int_{D_{\on{ref}}} a^{\mathrm{TT}}_{1, 1}(\hat{x}, \bs{y})\nabla^1\varphi_i(\hat{x})\nabla^1\varphi_{i'}(\hat{x})P_\alpha(\bs{y}) P_{\alpha'}(\bs{y}) \dx{\pi_M}(\bs{y})\dx{\hat{x}}.
\end{equation}
Using the multi-linear structure of $a^{\mathrm{TT}}_{1, 1}$, one can write $\mathbf{L}_1$ as
\begin{equation}
  \mathbf{L}_1(i, \alpha, i', \alpha') \approx \sum_{k_0=1}^{r_0^{1, 1}}\dots \sum_{k_{M-1}=1}^{r_{M-1}^{1, 1}} \mathbf{L}_0^1(i, i', k_0) \prod_{m=1}^M \mathbf{L}_m^1(k_{m-1}, \alpha, \alpha', k_m),
\end{equation}
where the first component tensor $\mathbf{L}^{1}_0$ depends on the physical discretization in piecewise constant FE functions $\{\psi_i\}_{i=1}^{N_0}$ only, i.e.,
\begin{equation}
 \mathbf{L}^1_0(i, i', k_0) = \sum_{\ell=1}^{N_0} a^0_{1, 1}(\ell, k_0)\int_{D_{\on{ref}}} \langle \nabla^1 \varphi_i(\hat{x}), \psi_\ell(\hat{x}) \nabla^1 \varphi_{i'}(\hat{x})\rangle \mathrm{d}\hat{x}.
\end{equation}
The remaining tensor operator parts decompose into one dimensional integrals over triple products of orthogonal polynomials of the form
\begin{equation}
  \mathbf{L}_m^1(k_{m-1}, \alpha, \alpha', k_m) = \sum_{\mu_m=0}^{d^{1, 1}_m-1} a^m_{1, 1}(k_{m-1}, \mu_m, k_m) \int_{[-1, 1]} P_{\mu_m} P_\alpha P_{\alpha'} \dx{\pi_m}.
\end{equation}
The evaluation is known explicitly thanks to the recursion formula for orthogonal polynomials, cf.~\cite{abra, ernst2012convergence}.

\begin{remark}
Due to the sum of TT operators in~\eqref{eq:TTop sum}, the result can be represented by a tensor with TT-rank: $d^2\max\{\mathbf{r}^{i, j}\;\vert\; i,j\in\{1, \ldots, d\}\}$.
\end{remark}


With the TT approximations of the data $f_{\mathrm{TT}} \approx \hat f$ and $\bs{a}_{i, j}\approx a_{i,j}^{\mathrm{TT}}$, we replace the original system of equations that have to be solved for $U\in\mathbb{R}^{N\times d_1\times \ldots\times d_M}$, namely
\begin{equation}
  \mathbf{L} U = \mathbf{F},
\end{equation}
with a constrained minimization problem on the low-rank manifold $\mathcal{M}_{\bs r}$ containing all tensor trains of dimensionality represented by $\Lambda$ and fixed rank $\bs r$,
\begin{equation}
  \label{eq:min problem}
  W^{\mathrm{TT}} = \argmin_{V\in\mathcal{M}_{\bs r}} \norm{\bs L^{\mathrm{TT}} V - \bs F^{\mathrm{TT}}}_{F}.
\end{equation}
Here, we take $\bs L^{\mathrm{TT}}$ and $\bs F^{\mathrm{TT}}$ as the TT approximations of $\bs L$ and $\bs F$, respectively and $\norm{\cdot}_{F}$ is the Frobenius norm.

To solve~\eqref{eq:min problem}, we chose a preconditioned alternating least squares (ALS) algorithm as described in~\cite{holtz2012alternating, EPS}.
This eventually results in an approximation of the Galerkin solution of~\eqref{eq:variational formulation} 
\begin{equation}
  \label{eq:approx solution}
  w_N \isdef w(\Lambda, \mathcal{T}, \bs r, \tau) = \sum_{\mu\in\Lambda}\sum_{i=1}^N W^{\mathrm{TT}}(i, \mu) \varphi_i P_\mu,
\end{equation}
where $\tau$ is a place-holder for the inscribed parameters of the numerical algorithm and $\bs r$ is the desired and predefined TT-rank of $W^{\mathrm{TT}}$.

%% file: adaptivity.tex
\section{Adaptive algorithm}
\label{sec:adaptive algorithm}

The error estimator of Section~\ref{sec:error estimator} is formulated in a computable TT representation in Section~\ref{sec:residual estimator}.
It gives rise to an adaptive algorithm, which refines the spatial discretization, the anisotropic stochastic polynomial set and the representation format iteratively based on local error estimators and indicators.
This enables the assessment of the development of the actual (unknown) error $\norm{u-w_N}_\mcA$.
The inherently high computational cost of the error estimators can be overcome by means of the tensor train formalism.
In what follows, we examine the efficient computation of the individual error estimator components in the TT format and describe the marking and refinement procedure.
For more details and a more general framework, we refer to the presentations in~\cite{EGSZ1, EGSZ2, EPS, eigel2018adaptive}.

\subsection{Efficient computation of error estimators}
\label{sec:residual estimator}

We illustrate the efficient computation of the deterministic error estimator $\eta_T$.
For each element $T\in\mathcal{T}$ of the triangulation, the error estimator is given by~\eqref{eq:est-cont}.
Due to the sum over $\Lambda$ it suffers from the curse of dimensionality.
However, employing the low-rank approximation $\bs A^{\on{TT}} \approx \bs A$, $f^{\on{TT}}\approx f$ and $w_N$ renders the computation feasible.
To make this more explicit, recall that
\begin{equation}
  \eta_T(w_N) = \norm{f^{\on{TT}} - \div(\bs A^{\on{TT}}\nabla_{\hat{x}}w_N)}_{L^2(\Gamma, \pi; L^2(T))}.
\end{equation}
This is evaluated by expansion of the inner product,
\begin{align}
  \label{eq:eta-sum}
  \eta_T(w_N)^2 = \norm{f^{\on{TT}}}^2_{L^2(\Gamma, \pi; L^2(T))} - &2\langle f^{\on{TT}}, \div(\bs A^{\on{TT}}\nabla_{\hat{x}}w_N)\rangle_{L^2(\Gamma, \pi; L^2(T))} \\ &+ \norm{\div(\bs A^{\on{TT}}\nabla_{\hat{x}}w_N)}^2_{L^2(\Gamma, \pi; L^2(T))}.
\end{align}
The first term is a simple inner product of a functional tensor train.
It reduces to a simple summation over the tensor components due to the orthonormality of the polynomial basis, i.e.,
\begin{equation}
  \norm{f^{\on{TT}}}^2_{L^2(\Gamma, \pi; L^2(T))} = \sum_{\mu\in\Lambda}\sum_{i'=1}^{N_0}\sum_{i=1}^{N_0}F(i, \mu)F(i', \mu)\int_T\psi_i\psi_{i'}\mathrm{d}\hat{x},
\end{equation}
whereas the high-dimensional sum can be evaluated for every tensor dimension in parallel using, for all $i, i'=1, \ldots, N_0$, that
\begin{align}
  \label{eq:compute F}
  \sum_{\mu\in\Lambda} F(i, \mu)F(i', \mu) = &\sum_{k_0=1}^{r^f_0}\ldots \sum_{k_{M-1}=1}^{r^f_{M-1}}\sum_{k_0'=1}^{r^f_0}\ldots \sum_{k_{M-1}'=1}^{r^f_{M-1}} F_0(i, k_0)F_0(i, k_0')\\ 
  &\prod_{m=1}^M \sum_{\mu_m=1}^{d_m}F_m(k_{m-1}, \mu_m, k_m)F_m(k_{m-1}', \mu_m, k_m').
\end{align}
Note that the iterated sum over the tensor ranks has to be interpreted as matrix-vector multiplications. 
Hence, the formula above can be evaluated highly efficiently.
In fact, if the employed TT format utilizes a component orthogonalization and $f^{\on{TT}}$ is left-orthogonal, the product can be neglected and one only has to sum over $k_0$ and $k_0'$. 

For the remaining terms in~\eqref{eq:eta-sum}, one has to find a suitable representation of $\bs A^{\on{TT}}\nabla_{\hat{x}}w_N$.
Since the gradient is a linear operator, one can calculate a tensor representation of this product explicitly, involving multiplied ranks and doubled polynomial degrees. 
For a detailed construction we refer to~\cite[Section 5]{eigel2018adaptive}.
The matrix-vector multiplication due to entry-wise TT representation of $\bs A^{\on{TT}}$ does not impose any further difficulties but a slight increase in complexity since one needs to cope with a sum of individual parts.
Eventually, the mixed and operator terms are computed in the same fashion, using the same arguments as for~\eqref{eq:compute F}.


\subsection{Fully adaptive algorithm}
\label{sec:fully adaptive algorithm}

Given an initial configuration consisting of a regular mesh $\mathcal{T}$, a finite active tensor multi-index set $\Lambda\subset \mathcal{F}$, a (possibly random) start tensor $W^{\mathrm{TT}}$ with TT-rank $\mathbf{r}$ and solver parameter $\tau$, consisting e.g. of a termination threshold, rounding parameter, iteration limit or precision arguments, we now present the adaptive refinement procedure summarized in Algorithm~\ref{alg:asgfem}.

On every level, we generate an approximation of the data $f^{\mathrm{TT}}$ and $\bs A^{\mathrm{TT}}$ by a tensor reconstruction. 
The procedure is e.g. described in~\cite{eigel2019non} and referred to as 
\begin{equation}
  f^{\mathrm{TT}},\; \bs A^{\mathrm{TT}} \leftarrow \Reconstruct[\Xi, \mathcal{T}, N_s],
\end{equation}
where the multi-index set $\Xi$ can be chosen arbitrarily, but it is advisable to consider Remark~\ref{rem:ttdim}. 
The number of samples $N_s$ can be related e.g. to Monte Carlo samples or more structured quadrature techniques such as Quasi Monte Carlo and sparse grid points.
In what follows we assume that the obtained approximations become sufficiently accurate.

The procedure for obtaining a numerical approximation $w_N\in\mathcal{V}_N$ is denoted by 
\begin{equation}
 w_N \leftarrow \Solve[\Lambda, \mathcal{T}, \tau, W^{\mathrm{TT}}].
\end{equation}
The used preconditioned ALS algorithm is only exemplary to obtain $w_N$. 
Alternative alternating methods or Riemannian algorithms are feasible as well.

To obtain the overall estimator $\Theta(\eta, \zeta,\iota)$, one has to evaluate the individual components by the following methods 
\begin{align*}
  (\eta_{T})_{T\in\mathcal T},\eta &\leftarrow \Estimate_x[w_N,f^{\mathrm{TT}}, \bs A^{\mathrm{TT}}, \Lambda,\mathcal T], \\
(\zeta_{m})_{m\in\mathbb N}, \zeta &\leftarrow \Estimate_{\bs y}[w_N,f^{\mathrm{TT}}, \bs A^{\mathrm{TT}}, \Lambda] \\
\iota &\leftarrow \Estimate_{\mathrm{LS}}[w_N, f^{\mathrm{TT}}, \bs A^{\mathrm{TT}}].
\end{align*}

A weighted balancing of the global estimator values $\eta, \zeta$ and $\iota$ results in the marking and refinement decision.

\subsubsection{Deterministic refinement}
In case of a dominant deterministic error estimator $\eta$, one employs a D\"orfler marking strategy on the mesh $\mathcal{T}$ for a ratio constant $\theta_\eta(w_N)$.
In abuse of notation, we use $(\eta_T)_{T\in\mathcal{T}}$ as the local error estimator on every triangle, where the jump components of $(\eta_S)_{S \in \mathcal{F}}$ are distributed among their nearby elements.
The method, consisting of the marking process and the conforming refinement of the marked triangles is covered by 
\begin{equation}
  \mathcal{T} \leftarrow \Refine_x[(\eta_T)_{T\in\mathcal{T}}, \eta, \mathcal{T}, \theta_\eta].
\end{equation}

\subsubsection{Stochastic refinement}
In case of a dominant stochastic error estimator $\zeta(w_N)$, we apply a D\"orfler marking on the set of local estimators $(\zeta_m)_{m\in\bbN}$ until the prescribed ratio $0 < \theta_\zeta < 1$ is reached. 
The marked dimensions in $\Lambda$ are increased $d_m \leftarrow d_m + 1$ by the method 
\begin{equation}
  \Lambda \leftarrow \Refine_{\bs y}[(\zeta_m)_{m\in\bbN}, \zeta, \Lambda, \theta_\zeta].
\end{equation}
\begin{remark}
 As stated in Section~\ref{sec:error estimator}, the global estimator $\zeta$ is not just the sum of the individual estimators $(\zeta_m)_{m\in\bbN}$ since the coupling structure is more involved. 
 Hence, we use $\zeta_{\mathrm{sum}}\isdef \sum_{m\in\bbN} \zeta_m$ in the marking procedure. 
 Due to the high regularity of the solution (Theorem~\ref{thm:regularity}), for $\Lambda$ large enough, one has $\zeta_{\mathrm{sum}} \approx \zeta$.
 Note that in the finite dimensional noise case, we have $\zeta_m = 0$ for $m>M$.
\end{remark}

\subsubsection{Representation refinement}
In the end, if $\iota$ has the largest contribution in the error estimator we improve the accuracy of the iterative solver.
For simplicity, we fix most of the solver parameter such as the number of alternating iteration or the termination value to low values that can be seen as overcautious. 
Nevertheless, in the low-rank tensor framework, the model class is restricted by the TT-rank $\mathbf{r}$.
Hence, we then allow $\mathbf{r} \leftarrow \mathbf{r} + \bs 1$ and add a random rank 1 tensor onto the solution tensor $W^{\mathrm{TT}}$. 
We summarize this approach in 
\begin{equation}
  W^{\mathrm{TT}}, \tau \leftarrow \Refine_{\mathrm{LS}}[W^{\mathrm{TT}}, \tau].
\end{equation}

\subsubsection{Adaptive algorithm}
One global iteration of this algorithm refines either the deterministic mesh $\mathcal{T}$, the active stochastic polynomial index-set $\Lambda$ or the tensor rank $\mathbf{r}$.
Iteration until the defined estimator $\Theta(\eta, \zeta,\iota)$ in Corollary~\ref{cor:error estimator} falls below a desired accuracy $\epsilon>0$ yields the adaptively constructed low-rank approximation $w_N\in\mathcal{V}_N$.

\begin{algorithm}
\begin{algorithmic}
\renewcommand{\algorithmicrequire}{\textbf{Input:}}
\renewcommand{\algorithmicensure}{\textbf{Output:}}
\REQUIRE{Initial guess $w_N$ with solution coefficient $W^{\mathrm{TT}}$; \\
solver accuracy $\tau$; \\
mesh $\mathcal T$ with degrees $p$; \\
active index set $\Lambda$;\\
sample size $N_s$ for reconstruction; \\
D\"orfler marking parameters $\theta_\eta$ and $\theta_\zeta$; \\
desired estimator $\Theta$ accuracy $\epsilon$.}
\ENSURE{New solution $w_N$ with new solution coefficient $W^+$; \\
\mbox{new mesh $\mathcal T^+$, or new index set $\Lambda^+$, or new tolerance $\tau^+$.}}
\REPEAT
\STATE{
\begin{tabular}{lll}
$f^{\mathrm{TT}}, \bs A^{\mathrm{TT}}$&$\leftarrow $&$\Reconstruct[\Xi, \mathcal T, N_s]$\\
$w_N $&$\leftarrow $&$\Solve[\Lambda,\mathcal T,\tau,W^{\mathrm{TT}}]$ \\
$(\eta_{T})_{T\in\mathcal T},\eta$ & $\leftarrow$ &$\Estimate_x[w_N,f^{\mathrm{TT}}, \bs A^{\mathrm{TT}}, \Lambda,\mathcal T,p]$ \\
$(\zeta_{m})_{m\in\mathbb N}, \zeta$ & $\leftarrow $&$ \Estimate_{\bs y}[w_N,f^{\mathrm{TT}}, \bs A^{\mathrm{TT}}, \Lambda]$ \\
$\iota $&$\leftarrow $&$\Estimate_{\mathrm{LS}}[w_N, f^{\mathrm{TT}}, \bs A^{\mathrm{TT}}]$
\end{tabular}}
\IF{$\max\{\eta,\zeta,\iota\} == \eta$}
  \STATE{$\mathcal T \qquad\;\, \leftarrow \Refine_x[(\eta_{T})_{T\in\mathcal T}, \eta, \mathcal T, \theta_\eta]$}
\ELSIF{$\max\{\eta,\zeta,\iota\} == \zeta$}
\STATE{$\Lambda \qquad\;\; \leftarrow \Refine_{\bs y}[(\zeta_{m})_{m\in\mathbb N}, \zeta, \Lambda, \theta_\zeta]$}
\ELSE 
\STATE{$W^{\mathrm{TT}},\tau \leftarrow \Refine_{\mathrm{LS}}[W^{\mathrm{TT}},\tau]$}
\ENDIF
\UNTIL{$\Theta(\eta, \zeta, \iota) < \epsilon$}
\RETURN{$w_N^+ = w_N;\; \mathcal{T}^+ = \mathcal{T};\; \Lambda^+ = \Lambda;\; \tau^+ = \tau$}
\caption{Reconstruction based adaptive stochastic Galerkin method}
\label{alg:asgfem}
\end{algorithmic}
\end{algorithm}

%% file: numericalResults.tex

\pgfplotstableread[col sep=comma]{
eps, r-est, r-adap, r-unif, e-m,               e-v,             m-dofs, tt-dofs, max-r
0.7, 0.46,  0.44,   0.45,   0.000190748367115, 0.0137347000359, 97346,  389384,   4
0.5, 0.46,  0.44,   0.45,   0.00383883107339,  0.0163227062458, 76740,  537285,   7
}\circle

\pgfplotstableread[col sep=comma]{
eps, r-est, r-adap, r-unif, e-m,                    e-v,                  m-dofs, tt-dofs, max-r
0.7, 0.51, 0.51,  0.46,   0.0070032354270719707,  0.021792384461861835, 81277,  406400,      5
0.5, 0.50, 0.51,  0.45,   0.0094526793869486384,  0.028698675709910054, 45050,  450760,     10
}\lshape

\pgfplotstableread[col sep=comma]{
eps, m1, m2
0.7, 2, 3
0.5, 5, 6
0.1, 21, 21
}\klmodes

\section{Numerical examples}
\label{sec:numerical examples}

This section is concerned with the demonstration of the performance of the described Galerkin tensor discretisation and the adaptive algorithm depicted in the preceding section.
We consider the linear second order model problem with a constant right-hand side and homogeneous Dirichlet boundary conditions

\begin{equation}\label{eq:modProbNumerics}
\begin{aligned}
-\div\big(\nabla u(\omega)\big) &= 1&&\text{in }D(\omega),\\
u(\omega) &= 0&&\text{on }\partial D(\omega),
\end{aligned}
\end{equation}
on two different reference domains in $\bbR^2$, namely the unit circle and the L-shape.
The Karhunen-Lo\`eve expansion of the random vector field stems from a Gaussian covariance kernel of the form
\begin{equation}
  \Cov[{\bs V}](\hat x, \hat{x}') = \frac{1}{1000}\begin{bmatrix}5\exp(-2\|\hat{x}-\hat{x}'\|^2_2) & \exp(-0.1\|2\hat{x}-\hat{x}'\|^2_2)\\   
\exp(-0.1\|\hat{x}-2\hat{x}'\|^2_2) & 5\exp(-0.5\|\hat{x}-\hat{x}'\|^2_2)\end{bmatrix}.
\end{equation}
The random variables in the Karhunen-Lo\`eve expansion are assumed to be independent and uniformly distributed on \([-\sqrt{3},\sqrt{3}]\), i.e.\ they
have normalized variance. Moreover, the mean is given by the identity, i.e.\ \(\E[{\bs V}](\hat{x})=\hat{x}\).

The computed spectral decomposition is truncated at a given threshold $\hat{\epsilon}$, which takes different values in the computational examples.
Table~\ref{tab:tolToM} summarizes how the choice of the truncation parameter affects the number of involved stochastic dimensions.
\begin{figure}
  \begin{center}
   \pgfplotstabletypeset[
   col sep=comma,
   every head row/.style={before row=\toprule,after row=\midrule},
   every last row/.style={after row=\bottomrule},
   columns/eps/.style={
     column name={tolerance $\hat{\epsilon}$},
   },
   columns/m1/.style={
     column name={circle KL terms},
   },
   columns/m2/.style={
     column name={L-shape KL terms},
   }
   ]{\klmodes}
  \caption{Comparison of different KL truncation tolerances and their
   implication on the number of KL modes $M$ for the employed reference domains.}
   \end{center}
     \label{tab:tolToM}
\end{figure}

We are interested in the correct approximation of the solution mean
\begin{equation}
  \mathbb{E}\left[u(\hat x, \cdot)\right] = \int_\Gamma u(\hat x, \bs y) \mathrm{d}\pi(\bs y)
\end{equation}
and solution variance 
\begin{equation}
  \mathbb{V}\left(u(\hat x, \cdot)\right) =   \mathbb{E}\left[u(\hat x, \cdot)^2\right] -   \mathbb{E}\left[u(\hat x, \cdot)\right]^2
\end{equation}
by means of the adaptive low-rank Galerkin approximation. 
In order to verify this, all experiments involve the computation of a reference mean and variance, based on a sampling approach. 
To that end, we employ the anisotropic sparse grid quadrature with Gauss-Legendre points\footnote{The implementation can be found online:  \texttt{https://github.com/muchip/SPQR}}, as described in \cite{HHPS18}.
The corresponding moments are then calculated on a fine reference mesh, resulting from uniform refinement of the last, adaptively computed, mesh $\mathcal{T}$, having at least $10^5$ degrees of freedom.
All experiments involve linear finite element spaces, i.e. $\mathcal{X}_p(\mathcal{T})$ with $p=1$. 
The number of quadrature points is chosen differently for the problems at hand. 
For $\hat{\epsilon}=0.7$ we take $53$ adaptively chosen nodes. 
Benchmarking the resulting mean from the sparse quadrature for this choice of samples against an approximation with additional nodes does not significantly improve the approximation quality (data not shown).
The same arguments apply for $\hat{\epsilon}=0.5$ and $301$ nodes, as well as for $\hat{\epsilon}=0.1$ and $4217$ nodes.
We denote the reference mean as $\mathbb{E}_{\on{ref}}[u]$ and the reference variance as $\mathbb{V}_{\on{ref}}(u)$.

\begin{remark}
  \label{rem:TTmean}
 In the low-rank tensor train format, the mean of a function, given in orthonormal polynomials, is computed highly efficiently, since the set of employed polynomials is orthonormal with respect to the constant function.
 Since, the corresponding coefficient is already incorporated in the representation, computing the mean is a simple tensor evaluation.
 More precisely, given $u\in\mathcal{V}_N$ we compute
\begin{equation}
   \mathbb{E}[u(\hat x, \cdot)] = \sum_{i=1}^N \sum_{\mu\in\Lambda} U(i, \mu) \varphi_i(\hat x) \mathbb{E}[P_\mu(\cdot)] = \sum_{i=1}^N U^{\mathrm{TT}}(i, \bs 0) \varphi_i(\hat x).
 \end{equation}
 Here, the evaluation of the tensor train $U^{\mathrm{TT}}$ at the multi-index $\bs 0=(0, \ldots, 0)$ consists of $M$ matrix-vector multiplications.
 
Similarly for the variance, we can compute the second moment as
\begin{equation}
   \mathbb{E}[u(\hat x, \cdot)^2] = \sum_{i=1}^N\sum_{i'=1}^{N} \sum_{\mu\in\Lambda} U(i, \mu) U(i', \mu) \varphi_i(\hat x)\varphi_{i'}(\hat x).
 \end{equation}
 This computation reduces even further, since in the tensor train setting of $U(i, \mu) = U_0(i)\otimes\bigotimes_{m=1}^M U_m(\mu_m)$, with $U_0\in\mathbb{R}^{N, r_0}$ and $U_m\in\mathbb{R}^{r_{m-1}, d_m, r_m}$ for $m=1,\ldots, M$ it is common to impose \emph{left-orthogonality}, i.e. $\sum_{\mu_m=1}^{d_m} U_m(\mu_m)U_m(\mu_m)^T = I_{r_{m-1}}$.
 Hence, the second moment reads
 \begin{equation}
   \mathbb{E}[u(\hat x, \cdot)^2] = \sum_{k_0=1}^{r_0}\sum_{i=1}^N\sum_{i'=1}^{N} U_0(i, k_0)U_0(i', k_0) \varphi_i(x)\varphi_{i'}(x),
 \end{equation}
 where it is advisable to not evaluate the matrix-matrix product over $k_0=1, \ldots, r_0$, since the resulting $\mathcal{O}(N^2)$ matrix is usually not sparse and exhibits available memory resources.
\end{remark}

The considered quantity of interest is the error of the mean and variance to the corresponding reference.
Therefore, for any TT approximation $w_N\in\mathcal{V}_N$ using remark~\ref{rem:TTmean}, we compute the relative error of the mean in $H_0^1(D_{\on{ref}})$ norm
\begin{equation}
  \label{eq:error quantity}
  e_E := \norm{\mathbb{E}_{\on{ref}}[u] - \mathbb{E}[w_N]}_{H_0^1(D_{\on{ref}})} \norm{\mathbb{E}_{\on{ref}}[u]}_{H_0^1(D_{\on{ref}})}^{-1}
\end{equation}
and the relative error of the variance in $W^{1, 1}(D_{\on{ref}})$ norm
\begin{equation}
  \label{eq:error quantity2}
  e_V := \norm{\mathbb{V}_{\on{ref}}(u) - \mathbb{V}(w_N)}_{W^{1,1}(D_{\on{ref}})} \norm{\mathbb{V}_{\on{ref}}(u)}_{W^{1, 1}(D_{\on{ref}})}^{-1}.
\end{equation}
Furthermore, we evaluate the rate of convergence of $e_E$ for the adaptive and uniform refinement, as well as the rate of the estimator $\Theta$. 
This is done by fitting the function $x\mapsto c_1e^{-\alpha}$ to the individual values, with respect to the number of degrees of freedom in the physical mesh. 
We denote the corresponding rates $\alpha_{E}^a$, $\alpha_E^u$ and $\alpha_\Theta$ respectively.

The employed tensor reconstruction algorithm is implemented in the open-source library \emph{xerus}~\cite{xerus}.
Every such approximated tensor is constructed on a set of $N_s$ samples $\{y^{(i)}\in\Gamma\}_{i=1\ldots N_s}$ as in the computation of the reference mean and polynomial degrees that are determined by the solution approximation as described in remark~\ref{rem:ttdim}.
In the considered examples, the tensor train solution employs constant and linear polynomials in all dimensions only. 
This behaviour is based on the stochastic estimator and reasoned in the complexity of the mesh approximation.
For the assembling of the physical part of the bilinear and linear form and the evaluation of sample solutions we make use of the PDE library \emph{FEniCS}~\cite{fenics}.
The entire adaptive stochastic Galerkin method is implemented in the open source framework \texttt{ALEA}~\cite{alea}.

\subsection{Example 1}
\label{sec:example1}
\begin{figure}
  \includegraphics[width=.3\linewidth]{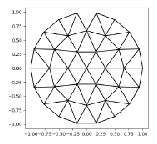}
  \includegraphics[width=.3\linewidth]{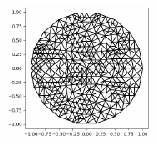}
  \includegraphics[width=.3\linewidth]{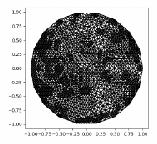}
  \caption{Adaptively refined mesh according to deterministic error estimator with mesh refinement ratio $\theta_\eta=0.2$ and KL truncation tolerance $\hat{\epsilon} = 0.1$.}
  \label{fig:mesh-circle}
\end{figure}

The first example is the random domain problem on the unit circle.
We use this problem as a reference, since the adaptive refinement is expected to yield similar results to uniform mesh refinement. 
Starting with an initial configuration of $16$ cells, fixed polynomial degree in the stochastic space of $d_1=\ldots=d_M=2$ and tensor rank $\mathbf{r}=\bs 2$, the described adaptive Galerkin FE algorithm yields the adaptively refined mesh depicted in Figure~\ref{fig:mesh-circle}.

\begin{figure}
\includegraphics[width=.49\linewidth]{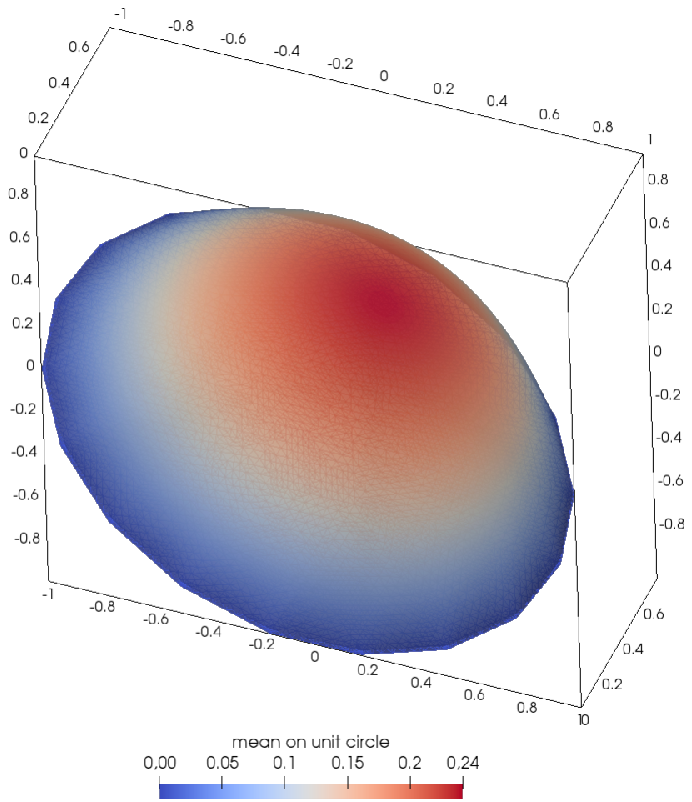}  
\includegraphics[width=.49\linewidth]{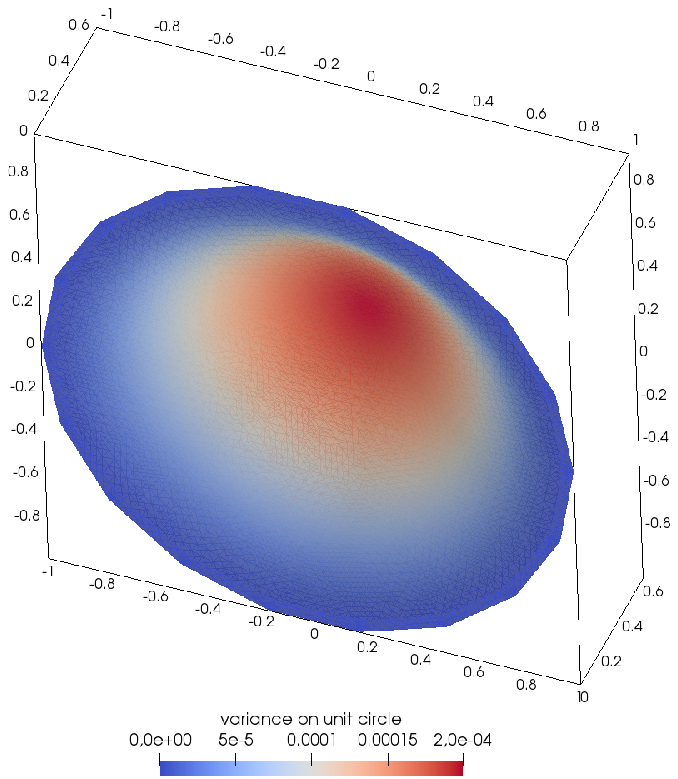}  
\caption{Computational mean (left) and variance (right) of the unit circle problem for KL tolerance $\hat{\epsilon}=0.1$.}
\label{fig:mean-circle}
\end{figure}

\begin{figure}
\centering
    \includegraphics[width=.2\linewidth]{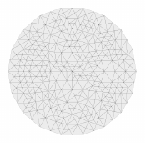}\\[-0.7em]
    
    \includegraphics[width=.2\linewidth]{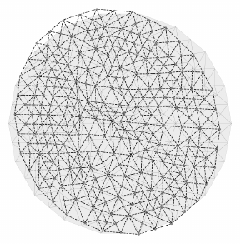}    
    \includegraphics[width=.2\linewidth]{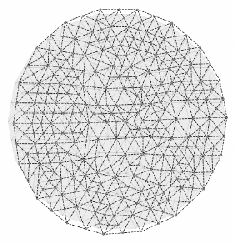}
    \includegraphics[width=.2\linewidth]{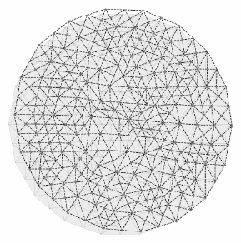}
    \includegraphics[width=.2\linewidth]{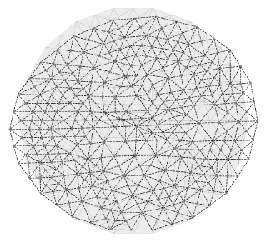}
    
    \includegraphics[width=.2\linewidth]{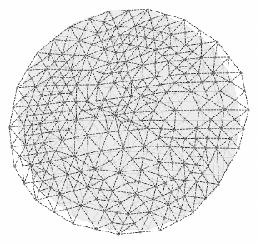}
    \includegraphics[width=.2\linewidth]{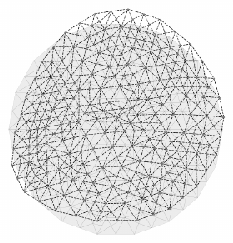}
    \includegraphics[width=.2\linewidth]{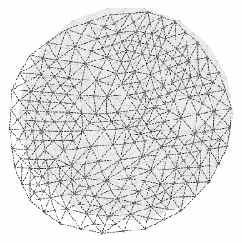}
    \includegraphics[width=.2\linewidth]{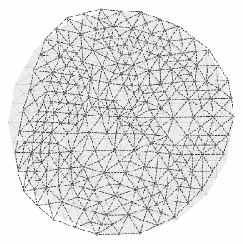}
\caption{Sample realisations of the transformed unit disc with KL tolerance $\hat{\epsilon}=0.5$ (middle-row) and $\hat{\epsilon}=0.1$ (bottom-row).
The reference mesh is given at the top.}
\label{fig:trafo-circle}
\end{figure}

For illustration purposes, we show the mean and variance of the solution on the unit disc together with realisations of the transformed reference domain for the adapted discretisation in Figures~\ref{fig:mean-circle} and~\ref{fig:trafo-circle}.
In Table~\ref{table:circle} we show the corresponding rates of convergence and the minimal reached error $e_E$ and $e_V$. 
The maximal involved tensor ranks are displayed in column $r_{\on{max}}$.
The degrees of freedom are shown for the physical mesh in column \emph{m-dofs} and for the tensor train itself in column \emph{tt-dofs}.
Note that the tensor train degrees of freedom refer to the dimension of the corresponding low-rank manifold.
As expected for the unit circle, the adaptive refinement does not improve the already optimal convergence rate. 

\begin{figure}
  \begin{center}
   \pgfplotstabletypeset[
   col sep=comma,
   every head row/.style={before row=\toprule,after row=\midrule},
   every last row/.style={after row=\bottomrule},
   columns/eps/.style={
     column name={$\hat{\epsilon}$},
   },
   columns/r-est/.style={
     column name={$\alpha_\Theta$},
   },
   columns/r-adap/.style={
     column name={$\alpha_E^a$},
   },
   columns/r-unif/.style={
     column name={$\alpha_E^u$},
   },
   columns/e-m/.style={
     column name={$e_E$},
   },
   columns/e-v/.style={
     column name={$e_V$},
   },
   columns/m-dofs/.style={
     column name={\emph{m-dofs}},
   },
   columns/tt-dofs/.style={
     column name={\emph{tt-dofs}},
   },
   columns/max-r/.style={
     column name={$r_{\on{max}}$},
   }
   ]{\circle}
   \caption{Results for the unit circle.
   Computed relative error of the mean and variance together with the rates of convergence for the error and the overall estimator $\Theta$. 
   Compared are the truncation tolerance $\hat{\epsilon}=0.7$ and $\hat{\epsilon}=0.5$.}
   \end{center}
   \label{table:circle}
\end{figure}

%
%
%
%
%

\subsection{Example 2}
\label{sec:example2}
For the second example, we chose the L-shaped domain $[-1, 1]^2\setminus\{[0, 1]\times [-1, 0]\}$.
The corner singularity is a typical example where adaptive refinement yields better approximation rates with respect to degrees of freedom than a uniform refinement.

\begin{figure}
  \includegraphics[width=.49\linewidth]{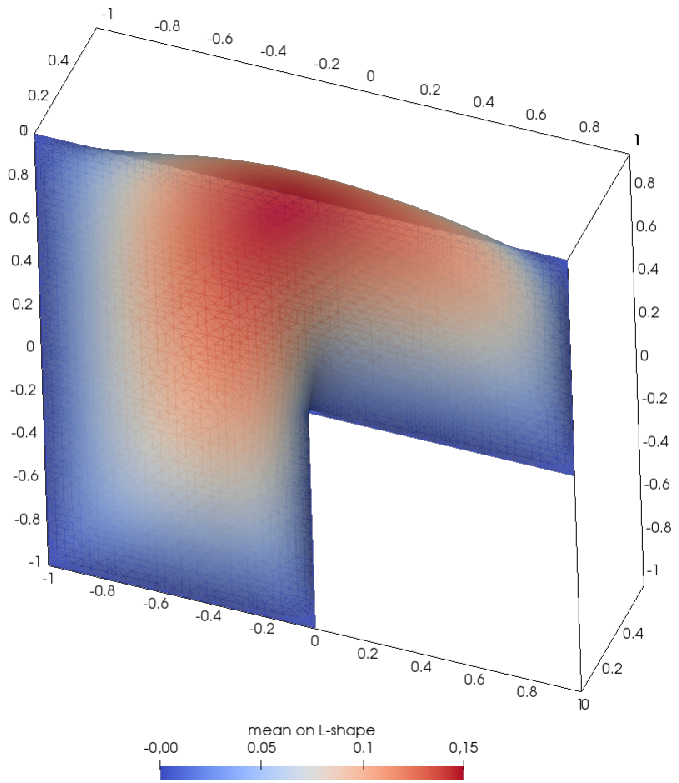}  
  \includegraphics[width=.49\linewidth]{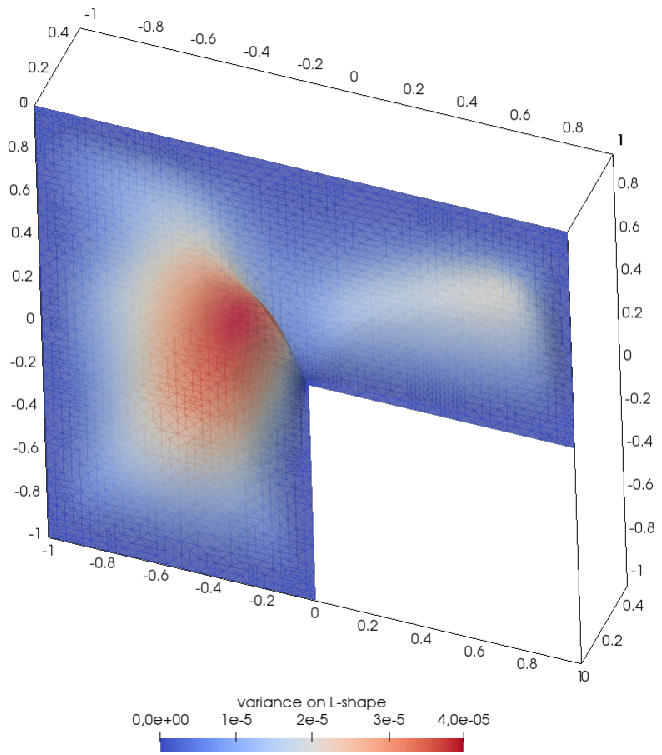}  
\caption{Computational mean (left) and variance (right) of the model problem on the L-shaped domain for KL tolerance $\hat{\epsilon}=0.5$.}
\label{fig:mean-lshape}
  
\end{figure}

\begin{figure}
  \includegraphics[width=.3\linewidth]{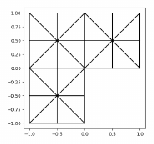}
  \includegraphics[width=.3\linewidth]{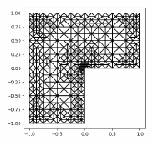}
  \includegraphics[width=.3\linewidth]{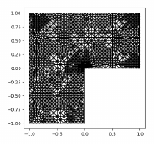}
  \caption{Adaptively refined mesh based on deterministic error estimator with mesh refinement ratio $\theta_\eta=0.5$ and KL truncation tolerance $\hat{\epsilon} = 0.5$.}
  \label{fig:mesh-lshape}
\end{figure}

\begin{figure}
\centering
    \includegraphics[width=.2\linewidth]{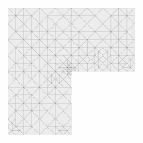}\\[-0.6em]
    
    \includegraphics[width=.2\linewidth]{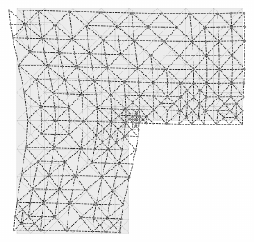}    
    \includegraphics[width=.2\linewidth]{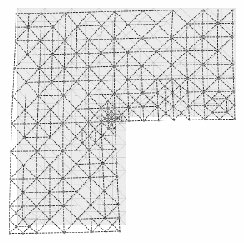}
    \includegraphics[width=.2\linewidth]{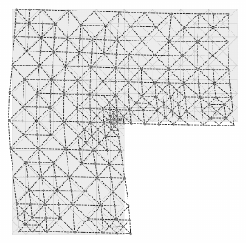}
    \includegraphics[width=.2\linewidth]{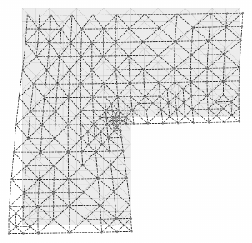}
    
    \includegraphics[width=.2\linewidth]{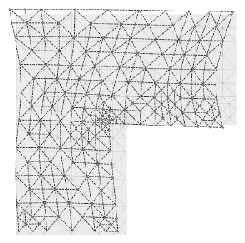}
    \includegraphics[width=.2\linewidth]{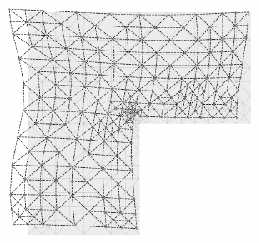}
    \includegraphics[width=.2\linewidth]{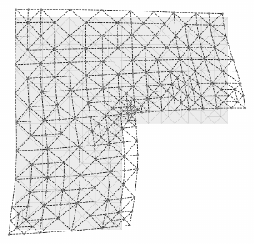}
    \includegraphics[width=.2\linewidth]{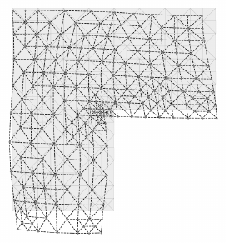}
\caption{Sample realisations of the transformed L-shape with KL tolerance $\hat{\epsilon}=0.5$ (middle-row) and $\hat{\epsilon}=0.1$ (bottom-row).
The reference mesh is given at the top.}
\label{fig:trafo-lshape}
\end{figure}

Starting with an initial configuration of $24$ cells, fixed polynomial degree in the stochastic space of $d_1=\ldots=d_M=2$ and tensor rank $\mathbf{r}=\bs 2$, the described adaptive Galerkin FE algorithm yields the adaptively refined mesh displayed in Figure~\ref{fig:mesh-lshape}.
Again, we show the approximated mean and variance together with some random realizations of the transformed reference domain in Figures~\ref{fig:mean-lshape} and~\ref{fig:trafo-lshape}.

%
%
%
%
%
%

\begin{figure}
\begin{center}
   \pgfplotstabletypeset[
   col sep=comma,
   every head row/.style={before row=\toprule,after row=\midrule},
   every last row/.style={after row=\bottomrule},
   columns/eps/.style={
     column name={$\hat{\epsilon}$},
   },
   columns/r-est/.style={
     column name={$\alpha_\Theta$},
   },
   columns/r-adap/.style={
     column name={$\alpha_E^a$},
   },
   columns/r-unif/.style={
     column name={$\alpha_E^u$},
   },
   columns/e-m/.style={
     column name={$e_E$},
   },
   columns/e-v/.style={
     column name={$e_V$},
   },
   columns/m-dofs/.style={
     column name={\emph{m-dofs}},
   },
   columns/tt-dofs/.style={
     column name={\emph{tt-dofs}},
   },
   columns/max-r/.style={
     column name={$r_{\on{max}}$},
   }
   ]{\lshape}
   \caption{Results for the L-shape domain. 
   Computed relative error of the mean and variance together with the rates of convergence for the error and the overall estimator $\Theta$. 
   Compared are the truncation tolerance $\hat{\epsilon}=0.7$ and $\hat{\epsilon}=0.5$.}
   \end{center}
   \label{table:lshape}
\end{figure}  

The obtained rates of convergence and error quantities are shown in Table~\ref{table:lshape}.
Fortunately, the obtained rate for the estimator $\alpha_\Theta$ follows the error rate $\alpha_E^a$, in contrast to the uniform refinement strategy, where a slightly slower convergence is achieved.